\documentclass[14pt,reqno,twoside]{amsart}
\usepackage{amsmath}
\usepackage{amssymb}
\usepackage{graphicx}

\newcommand{\ignore}[1] {}

\newcommand{\R}{{\mathbb{R}}}
\newcommand{\Z}{{\mathbb{Z}}}

\newcommand{\N}{{\mathbb{N}}}

\newcommand\NN{\mathbb N}
\newcommand\RR{\mathbb R}
\newcommand\ZZ{\mathbb Z}
\newcommand\QQ{\mathbb Q}

\newcommand\cA{\mathcal{A}}
\newcommand\cB{\mathcal{B}}

\newcommand\cD{\mathcal{D}}

\newcommand\cM{\mathcal{M}}
\newcommand\cS{\mathcal{S}}

\newcommand\cU{\mathcal{U}}

\newcommand\cL{\mathcal{L}}

\newcommand\cG{\mathcal{G}}

\newcommand\norm[1]{\left\|#1\right\|}

\newcommand\abs[1]{\left|#1\right|}

\newcommand\inn[1]{\left\langle #1 \right\rangle}
\newcommand\set[1]{\left\{{#1}\right\}}

\newtheorem{thm}{Theorem}
\newtheorem{lem}[thm]{Lemma}
\newtheorem{prop}[thm]{Proposition}
\newtheorem{cor}[thm]{Corollary}
\newtheorem{remark}[thm]{Remark}
\newtheorem{definition}[thm]{Definition}
\newtheorem{example}[thm]{Example}
\newtheorem{question}[thm]{Question}

\begin{document}

\title[Hardy-Littlewood maximal inequality]{Hardy-Littlewood maximal 
operator on spaces with exponential volume growth}
%

\author{Koji Fujiwara}
\address{Department of Mathematics, Kyoto University}
\email{kfujiwara@math.kyoto-u.ac.jp}
\thanks{The first author was supported in part by Grant-in-Aid for Scientific Research (20H00114)}

\author{Amos Nevo}
\address{Department of Mathematics, University of Chicago, and Technion, IIT}
\email{amosnevo6@gmail.com}


\date{\today}



\begin{abstract}

We consider the Hardy-Littlewood maximal function associated with ball averages on  spaces with exponential volume growth. We focus on discrete groups with balls defined by invariant metrics associated with a variety of length functions. Under natural assumptions on the rough radial structure of the group in question, we establish a weak-type $\cL\left(\log \cL\right)^{\bf c}$ maximal inequality for the Hardy-Littlewood maximal function.  We give a variety of examples where the rough radial structure assumptions hold, based on considerations from geometric group theory, or   on analytic considerations related to the regular representation of the group. We elucidate the connections of these assumptions to a spherical coarse median inequality, to almost exact polynomial-exponential growth of balls, and to the radial rapid decay property. 

In particular, the weak-type maximal inequality in $\cL\left(\log \cL\right)^{\bf c}$ is established for any lattice in a connected semisimple Lie group with finite center, with respect to the distance function restricted from the Riemannian distance on symmetric space to an orbit of the lattice. It is also established for right-angled Artin groups, Coxeter groups and braid groups, for a suitable choice of word metric. 

For non-elementary word-hyperbolic group we establish that the Hardy-Littlewood maximal operator with respect to balls defined by a word length satisfies the weak-type $(1,1)$ maximal inequality, which is the optimal result. 
\end{abstract}

\maketitle
\section{Ball averages and maximal inequalities}
\subsection{Introduction} 
Let $X$ be a locally-compact second-countable (lcsc) non compact space. Let $D_X$ be a semi-metric on $X$ such that the balls given by  $B_r(x_0)=\set{x\in X\,;\, D_X(x,x_0) <  r}$ are pre-compact, for al $x_0\in X$. Let $m_X$ be an unbounded Radon measure on $X$, so that the measure of every compact subset is finite. We assume that for some fixed positive constant $r_0$, the measure of the balls $B_r(x_0)$, $ r > r_0 $ is at least a fixed positive constant, uniformly for all $x_0 \in X$.  For every radius $r > r_0$ define the ball averaging operator on the space of locally integrable functions $f$ on $X$ via 
\begin{equation}\label{def-ball-ave} 
 \cB_r(f)(x_0)=\frac{1}{m_X(B_r(x_0))}\int_{x\in B_r(x_0)} f(x)dm_X(x)
\end{equation} 
Define the maximal operator associated with the family of ball averages with $r > r_0$ by 
\begin{equation}\label{def-HL-max-op}
\cM f(x_0)=\sup_{r  > r_0}  \cB_r\left(\abs{f}\right)(x_0)
\end{equation}
The weak-type $(1,1)$-maximal inequality for the (global part of the) Hardy-Littlewood maximal operator $\cM$ asserts that there exists a finite constant $C$, independent of $f$, such that for every $\eta > 0$
  \begin{equation}\label{def-HL-max}
m_X\set{x\,;\, \cM f(x)\ge \eta}\le \frac{C}{\eta}\norm{f}_{L^1(X)}\,.
\end{equation}
We can now formulate :

{\bf The Hardy-Littlewood problem for ball averages.}
When does the operator $\cM(f)$ satisfy the weak-type $(1,1)$-maximal inequality ?

To develop a productive approach to this problem, it seems necessary to start by assuming, at the very least, the following two assumptions. First, the two families of balls $B_r(x_0)$ and $B_r(x_0^\prime)$ with $x_0, x_0^\prime$ a bounded distance apart are uniformly comparable in shape in a suitable sense, one such being that they are almost isometric to one another. Second, the families of balls satisfy certain conditions on their volume growth, in a uniform fashion. Our main focus below will be on the Hardy-Littlewood maximal inequalities of weak-type  on spaces of exponential growth, and for a discussion of certain such spaces,  under the assumption of bounded geometry, we refer to  \cite{MPSV} (for manifolds), \cite{LMSV} (for trees), and \cite{LZ23} (for products of hyperbolic spaces with manifolds with the volume doubling property).

Rather than assuming such conditions directly, we will consider natural  geometric structures invariant under a group action, that will guarantee the uniformity properties mentioned above. Two sets of assumptions based on the existence of geometric symmetries of the space are of particular significance and are given as follows. 
\begin{itemize} 
\item $X=G/L$ is a homogeneous space under a transitive action of a locally-compact second-countable (lcsc) group, and furthermore $G$ acts isometrically w.r.t. to the metric $d_X$, and preserves the measure $m_X$, defined using a left-invariant Haar measure on $G$. 
\item  $\Gamma$ is a countable group which acts properly discontinuously and co-compactly on $X$, isometrically w.r.t. to the metric $D_X$, and preserves the measure $m_X$. 
\end{itemize}

 To gauge the extensive scope of the ball averaging problems under these assumptions, note that  it includes the following.
 \begin{enumerate}
\item Homogeneous spaces $G/K$ for all connected Lie groups $G$ and their compact subgroups $K$, including for example $K=\set{e}$, with the measure $m_X$ defined by a left-invariant Riemannian structure, and the metric defined either by the Riemannian metric or by a length function associated to a $K$-invariant norm in a linear representation of $G$ with compact kernel. 
 \item Homogeneous spaces, and more generally spaces with finitely many orbits,  arising from totally disconnected groups and their compact subgroups. This includes, in particular Bruhat-Tits buildings associated with semisimple algebraic groups over local fields. 
 \item The universal Riemannian covering space $X$ of a closed Riemannian manifold $M$, with the fundamental group $\Gamma=\pi_1(M)$ acting isometrically, preserving the Riemannian metric and the Riemannian volume.  
  \item Finitely-generated countable groups $\Gamma$, acting by left translations on $X=\Gamma$, isometrically with respect to an invariant metric. Such metrics can be given by a symmetric word metric associated with finite generating set,  or by a length function defined by a norm in a faithful unbounded linear representation of $\Gamma$. The invariant measure is the counting measure.  
\end{enumerate}

{\it The local Hardy-Littlewood operator}. In the examples just noted, the local part of the Hardy-Littlewood maximal operator, given by the supremum of $\cB_r(\abs{f})$ over $0 < r \le r_0$, can easily be analyzed. When $X$ is a homogeneous space of the groups mentioned above, or a Riemanian manifold, the local averaging operator $\cB_r$, $0<r \le r_0$ can be estimated using the (local) doubling condition for balls with bounded radius which holds in these cases, and the general local maximal inequality can be established using a local version of the transference principle.
We refer to  \cite[\S 6.1, Thm. 6.2]{N06} for a full discussion.  

\subsection{Evolution of the problem} 
Let us begin by very briefly noting some of the history of the Hardy-Littlewood maximal inequality and some of the main results obtained. 

First, given the invariant (semi-)metric $D_X$ and the invariant Radon measure 
$m_X$, recall the two scales measuring the volume growth of balls centered at $x\in X$, as follows.

The polynomial scale is defined by : 
\begin{equation}\label{poly-growth} 
\mathfrak{p}(X,x)=\limsup_{r\to \infty}\frac{ \log m_X(B_R(x))}{\log r}\,.
\end{equation} 
When $\mathfrak{p}(X,x)$ has a global finite upper bound which is independent of $x\in X$ the system $(X, D_X, m_X)$ is called of polynomial growth. The least upper bound 
of $\mathfrak{p}(X,x)$ for $x\in X$ is then called the degree of polynomial volume growth, denoted $\mathfrak{p}(X)$.

The exponential scale is defined by: 
\begin{equation}\label{vol-growth}
\mathfrak{g}(X,x)=\lim_{r\to \infty}\frac{1}{r} \log m_X(B_r(x)).
\end{equation} 
Under our assumptions, in all the cases we consider below, the limit does in fact exist and the growth parameter is independent of $x\in X$, so denote it by $\mathfrak{g}=\mathfrak{g}(X)$. 
When $\mathfrak{g}(X) > 0$ the system $( X, D_X, m_X)$ is called a system of exponential volume growth. When $\mathfrak{g}=0$ it is called of subexponential growth.

\subsection*{I. Systems with polynomial volume growth}
The origin of the weak-type $(1,1)$-maximal inequalities for ball averages is the  classical theorem of Hardy and Littlewood \cite{HL30}, establishing it for the real line, with respect to the usual metric and Lebesgue measure. Averages with respect to Euclidean balls in Euclidean spaces of arbitrary dimension were considered by Wiener \cite{Wi39}, who established the weak-type $(1,1)$ maximal inequality by means of the Wiener covering argument. 

A crucial fact underlying the covering argument for Euclidean balls is the volume doubling property discussed in general by Calderon \cite{C53}, namely the condition that 
$m_X(B_{2r})\le C m_X(B_r)$ (see the discussion in \cite[Ch.1]{St93}). This  condition is certainly satisfied whenever the volume of the balls satisfies  $C_1 r^d \le m_X(B_r)\le C_2  r^d$, and so in particular holds for a large collection of metrics on connected nilpotent  groups $X=G$, where the volume satisfy this growth condition (see \cite{Gu73} for this estimate, and  \cite{Pa83} and \cite{Br05} for the existence of the limit  $m_X(B_r)/r^d)$).  The class of lcsc groups $G$ (and homogeneous spaces $G/K$  with $K$ compact) with polynomial volume growth is the only general class for which the Hardy-Littlewood weak-type $(1,1)$-maximal inequality has been established (for an account, see  \cite[Thm. 5.7]{N06} and the discussion in \cite[\S 5]{N06}). For an extensive discussion of applications of the weak-type maximal inequality for ball averages in this context and its significance in real-variable theory, together with many generalizations and ramifications,  we refer to \cite[Ch. I]{St93}. 

\subsection*{II. Systems with exponential volume growth.}

When the space $X$ has exponential volume growth, the Hardy-Littlewood weak-type maximal inequality presents a notoriously difficult analytic problem. Not much is known  about it, and the main results in the specific geometric group-invariant context we are considering that we are aware of are as follows. 

{\bf 1) Symmetric spaces.} Let $G$ be a connected non-compact semisimple Lie group with finite center, and $K$ a maximal compact subgroup. $X=G/K$ has the structure of a symmetric space, with $G$ acting transitively and isometrically on $X$, preserving a Riemannian volume. Examples include hyperbolic spaces (associated with simple groups of real-rank one), as well  symmetric spaces of non-positive curvature, associated with semi-simple Lie groups with real-rank at least $2$, such as $SL_n(\RR)$, $n \ge 3$. 
The Hardy-Littlewood weak-type $(1,1)$-maximal inequality for the ball averages defined using the Riemannian metric associated with the Killing form was established for all symmetric spaces by Str\"omberg \cite{Str81}. The method proceeds by passage to the Iwasawa horospherical coordinate system $G=KAN=KP$, where $P$ is the stability group of a point in the maximal boundary $G/P$. This coordinate system is used to parametrize radial averages on the symmetric space, rather than the obvious radial coordinates, and constitutes a crucial tool in the proof.  

{\bf 2). Exponential solvable groups  $S$ :  Harmonic AN groups.} Connected simply-connected exponential solvable groups $S$ have exponential volume growth, and the Hardy-Littlewood ball averaging problem was considered for a special class of them - the harmonic $AN$ groups. This subclass generalizes the groups $AN$ appearing in the Iwasawa decomposition of connected simple Lie groups $G=KAN$ of real rank one and finite center, and exhibit remarkably similar geometric properties.  In particular, radial harmonic analysis on the symmetric space $X$ can be generalized to  harmonic $AN$ groups, and in \cite{ADY96} Anker, Damek and Yacoub developed this analogy in full, and adapted St\"omberg's approach to prove the weak-type $(1,1)$ maximal inequality for Riemannian ball averages on them (see \cite[Thm. 3.14, Cor. 3.22]{ADY96}).

{\bf 3). The homogeneous tree and the free group.} 
Consider the homogeneous tree $T$ with valency at least $3$, and with the usual graph metric and the counting measure. When the valency is even, say $2k$, there is a free group on $k$ free generators acting simply-transitively on its vertices, and the tree can be equivariantly identified with the Cayley graph of the free group with respect to a set of $k$ free generators. In this set-up the Hardy-Littlewood weak-type $(1,1)$ maximal inequality for the ball averages has been given several proofs, as follows.  
\begin{enumerate}
\item  The first proof was given by Rochberg and Taibleson in 1991 \cite{RT91}, using random walk methods and explicit computation of the Green function for the uniform random walk on the regular tree. Such an explicit computation of the Green function  is also possible in a few other cases, including the semi-homogeneous tree, and the same method then applies. This implies that several systems similar to the homogeneous tree and the free groups also satisfy the maximal inequality. Such systems include certain Cayley graphs of certain free products $G_1\ast G_2$ of two finite groups or $H_1\ast H_2\ast\cdots \ast H_N$ where $H_i$, $1\le i\le N$ are finite groups of equal size. The groups in question act on the semi-homogeneous tree or on certain distance-transitive graphs properly and co-finitely (see \cite{Mac82} for the connection between these two structures). 
\item Another proof for the homogeneous tree was given by Cowling, Meda and Setti and eventually published in 2009 \cite{CMS}.
Here the proof uses the identification of the homogeneous tree as a homogeneous space of its isometry group $\text{Isom}(T)$ modulo the compact stability group $K$ of a vertex, as well as the fact that its boundary can likewise be identified with the homogeneous space $\text{Isom}(T)/P$, where $P$ is the stability group of a boundary point. 
The horospherical coordinates on the tree are used along the lines of the arguments of 
Str\"omberg and Anker-Damek-Yacoub in the real-hyperbolic and harmonic-$AN$ cases, respectively. Again, this proof applies somewhat more generally, for example to the case of the semi-homogeneous trees and  to certain lattice subgroups of its isometry groups. 
\item  The rooted (and the homogeneous) $k$-tree were analyzed by Naor-Tao in 2009 \cite{NT09}, who gave another proof of the weak-type $(1,1)$ maximal inequality. This proof constitutes an important new departure, since it develops several novel geometric estimates based solely on the {\it radial structure} of the homogeneous tree. These arguments were the starting point of the present study, as we will point out in detail below.  
We remark that other applications of the methods developed  \cite{NT09} were considered in \cite{ORRS} and \cite{ST19}. 
\end{enumerate}

 \subsection*{III. Systems with intermediate volume growth.} 
As is well known, some countable finitely-generated groups $\Gamma$ have intermediate word growth : when $\Gamma$ acts on $X=\Gamma$  by left translations w.r.t. a symmetric word metric and the counting measure, the balls have 
both sub-exponential growth  namely $\mathfrak{g}(X)=0$, and also super-polynomial growth namely 
$\mathfrak{p}(X)=\infty$. 

We are not aware of any results on the Hardy-Littlewood maximal operator in this set-up.

\subsection{Rough growth and structure conditions}
In the present paper we focus on  countable finitely generated groups, acting by left translation, isometrically with respect to word metrics or other (suitably chosen) length functions $\cG: \Gamma \to \RR_{\ge 0}$. We will focus exclusively on groups with exponential growth, and anticipating our arguments below, we begin by defining spherical shells (sometimes also referred to as annuli) in the group. Namely we fix $L> 0$ and consider the sets 
\begin{equation} 
SS_r=\set{\gamma\in \Gamma\,;\, r\le \cG(\gamma)< r+L}.
\end{equation} 
For simplicity of notation, we will suppress the explicit dependence of the family of spherical shells on the parameter $L$.

We will present a method of establishing weak-type maximal inequalities for an extensive class of groups with exponential growth, based only on (suitable combinations of) the following general "radial rough structure" assumptions (which will be discussed in detail below).

\begin{enumerate}
\item {\bf Almost exact growth condition :} The growth of the size of the spherical shells with respect to the invariant metric is almost exactly polynomial-exponential, namely bounded above and below by constant multiples of $ r^{d} q^r$, for some $d \ge 0$, $q> 1$, when $r \ge 1$. 
\item {\bf  Rapid decay of averages on spherical shells: }The averaging operators on spherical shells, acting in $\ell^2(\Gamma)$, satisfy the spectral condition of rapid decay of their convolution norms. 
\item {\bf Rough structure condition:} The group satisfies a spherical coarse-median inequality (defined in eq. (\ref{median-structure}) below).  
\end{enumerate}

We note that these conditions are formulated solely in terms of the underlying invariant distance and the size of the spherical shells.  Our proof will proceed by establishing a weak-type maximal inequality for averages on spherical shells, and then deducing a weak-type maximal inequality for the ball averages.

We note that the fact that it is possible to use the spherical shells averages to control the ball averages is a reflection of the fact that the groups we consider have exponential volume growth, and more precisely of the fact that the size of $SS_r$ is at least a fixed positive fraction of the size of the ball $B_r$. Note that in contrast,  an approach attempting to deduce  the maximal inequality for ball averages from a maximal inequality for spherical shell averages is bound to fail when the volume growth is polynomial. Indeed in this case the size of the spherical shells becomes negligible compared with the size of the balls as the radius grows, and the family of spherical shell averages does not satisfy the same maximal inequalities as those satisfied by the family of ball averages (as exemplified, for instance, by the case of the groups $\ZZ^d$, see \cite{N04} for a full discussion).

Our main results  will show that all groups which satisfy the growth condition (1) and either one of conditions (2) and (3) stated above, with respect to a suitable length function (including all symmetric word metrics), satisfy a  weak-type maximal inequality for ball averages in the function space $\cL\left(\log\cL\right)^{\bf c}$ defined below. Here  the polynomial parameter $\bf c$ will be explicitly given in terms of some parameters appearing in the rough structure assumptions. Under conditions (1) and (3) we will establish in some cases that $\bf c=0$ namely that the weak-type $(1,1)$ maximal inequality holds for the Hardy-Littlewood operator.
  
We will also give an extensive list of examples of classes of groups for which the key assumptions above are satisfied, see Section \ref{sec.example}. 

We now turn to define and discuss our rough-structure assumptions. 

\subsection{Groups of polynomial-exponential  growth}

Let $\Gamma$ be a countable group, and let $\cG: \Gamma \to \RR_{\ge 0}$ be a   length function as defined in the introduction, namely a proper function satisfying $\cG(uv)\le \cG(u)+\cG(v)$, which is symmetric, namely $\cG(u^{-1})=\cG(u)$, and satisfied $\cG(e)=0$. 

Then $\cG$ defines the left-invariant semi-metric associated to the length function  by $d_\cG(u,v)=\cG(u^{-1} v)$, which is a metric if and only if the length function is faithful, namely vanishes at $g=e$ only.  
The associated balls centered at $e$ are given by
$$B_r=\set{\gamma\in \Gamma \,;\, \cG(\gamma)<  r}\,.$$ 
Each ball is a finite set, since the length function is a proper function. 
Spherical shells of width $L$ are given by 
$$SS_r=\set{\gamma\in \Gamma \,;\, r\le \cG(\gamma)< r+L}\,.$$ 
Here $L$ is a fixed positive constant, which will be chosen below to be suitably large.

  For simplicity we suppress it from the notation. In the special case of an integer-valued gauge, the spheres $S_r=\set{w\in \Gamma \,;\, \cG(w)=r}$ (where $r \ge 0$) coincide with the spherical shells $SS_r$ for any choice $0< L\le 1$.

\begin{definition}[]\label{exact-growth}
\begin{enumerate}
\item 
$(\Gamma, \cG)$ has {\em almost-exact polynomial-exponential growth} 
 if there exists $q > 1$, $0\le  d\in \NN$, $L_0 > 0$ and $C_{gr} \ge 1$ such that for all $r\ge 1$ 
 the size of the family of spherical shell $SS_r$  (for any fixed width $L$ with $L\ge L_0 > 0$) satisfies 
\begin{equation}\label{exact growth} 
C_{gr}^{-1} r^{d} q^r\le \abs{SS_r}\le C_{gr} r^{d} q^r\,.
\end{equation}
Here $C_{gr}$ depends of the width $L$ of the spherical shells, but the parameters $(q,d)$, when they exist, are uniquely determined and are referred to as the polynomial-eponential  growth parameters of the length function $\cG$.

If $d=0$, then $(\Gamma, \cG)$ has {\em almost exact exponential growth}.  
\item $(\Gamma, \cG)$ has  {\em exact polynomial-exponential growth} if
$$\lim_{r \to \infty} |SS_r|/(r^{d} q^r)$$ exists and is finite and positive. 
\end{enumerate} 
\end{definition}
\begin{remark}{\bf A note about terminology}. We have used the term "almost exact" growth in order to evoke the concept of "almost isometry" between metric spaces. This is a far more stringent condition than "quasi isometry", and the one suitable to discussing the quantitative aspects of growth in groups. Indeed, let $(\Gamma_1, \cG_1)$ and $(\Gamma_2, \cG_2)$ be almost isometric via a map $\Phi$, namely 
$$d_1(\gamma_1,\gamma_2)-C \le d_2(\Phi(\gamma_1),\Phi(\gamma_2))\le d_1(\gamma_1,\gamma_2)+C. $$
 Then if one of the spaces has almost exact polynomial-exponential growth, then so does the other, and the growth parameters $(d,q)$ are equal. 
We note that "almost exact growth" was called "strict growth" in \cite[\S 4. Def. 4.3]{N06}. 
\end{remark}

\subsection{Radial averages as convolutions}

Given any two functions $f,h \in \ell^1(\Gamma)$ their convolution is the function 
$f\ast h\in\ell^1(\Gamma)$ defined by 

\begin{equation}\label{conv}
f\ast h(w)=\sum_{u,v\in \Gamma\,;\, uv=w} f(u)h(v)
\end{equation}
$$=\sum_{u\in \Gamma} f(u)h(u^{-1}w)
=\sum_{v\in \Gamma} f(wv^{-1})h(v). $$
Let $A \subset \Gamma$ be a non-empty finite set. Clearly 
$\norm{1_A}_{\ell^1}=\abs{A}$ and $\norm{1_A}_{\ell^2}=\sqrt{\abs{A}}$.
Given another non-empty finite set $B$ we have 
\begin{equation}\label{conv-estim-1}
1_A \ast 1_B(w)=\sum_{u,v\in \Gamma\,;\, uv=w} 1_A(u)1_B(v)=
\abs{\set{u\in A, v\in B\,;\, uv=w}}
\end{equation}

so that  
\begin{equation}\label{conv-estim-2}
\norm{1_A \ast 1_B}^2_{\ell^2}=\sum_{w\in \Gamma} \abs{1_A \ast 1_B(w)}^2= 
\sum_{w\in \Gamma}\abs{\set{u\in A, v\in B\,;\, uv=w}}^2.
\end{equation}

The operator of averaging the function $h$ over a (non-empty) sphere of radius $r\ge 0$ on $\Gamma$ is given by (upon setting $u=wv$ and noting that $S_r$ is symmetric) 
\begin{equation}\label{radial conv} \cS_r h(w)=\frac{1}{\abs{S_r} } \sum_{u\,;\, d_\cG(w,u)=r}  h(u)=\frac{1}{\abs{S_r} } \sum_{u\,;\, \cG(w^{-1} u)=r}  h(u)
\end{equation}
$$=\frac{1}{\abs{S_r} } \sum_{v\in S_r}  h(wv)\sigma_r(v^{-1})=h\ast \sigma_r(w),$$
 Therefore $\cS_r$ coincides with the operator of right-convolution by the probability measure $\sigma_r$, the normalized sphere average. Every radial average is defined by a probability measure on $\Gamma$ which a convex combinations of the averages $\sigma_r$, and again is given by a right convolution operator.  We denote the normalized ball average by $\beta_r$ and the corresponding operator on $\ell^2(\Gamma)$ by $\cB_r$, and the normalized spherical shell average is denoted by $\sigma\sigma_r$, with the corresponding operator on $\ell^2(\Gamma)$ denoted by $\cS\cS_r$.

\section{Radial rapid decay and convolution estimates} 
Let us recall the following well-known definitions introduced in \cite{J}. We refer to e.g.  \cite{CR05} and \cite{Ch17} for a fuller discussion. 

\subsection{Rapid decay and radial rapid decay}

Consider the group $\Gamma$ with the length function 
$\cG$, and let $f$ be a given finitely supported function in $\Gamma$. For all $h\in \ell^2(\Gamma)$  
consider the inequality (where ${\bf b}\ge 0,C_{rd}>0$ are two fixed constants, depending only on $\Gamma$ and $\cG$)
\begin{equation}\label{RD-def-1}
\norm{f\ast h}_{\ell^2}\le C_{rd}  \, r^{\bf b} \norm{f}_{\ell^2} \norm{h}_{\ell^2}.
\end{equation} 
Equivalently, if $f\ast$ denotes the convolution operator defined by $f$, consider the inequality: 
\begin{equation}\label{norm-estimate-RD}
\norm{f\ast }_{\ell^2\to \ell^2}\le C_{rd}  \, r^{\bf b} \norm{f}_{\ell^2} .
\end{equation} 

\begin{definition}\label{RD-def} 
 Let $r\ge 1$. 
\begin{enumerate}
\item If for any given function $f$ supported in the ball $B_r(e)$  and every $h\in \ell^2(\Gamma)$  (\ref{RD-def-1}) is satisfied, then $(\Gamma,\cG)$ satisfies the property of {\bf rapid decay} with polynomial parameter ${\bf b}$. 
\item  If (\ref{RD-def-1}) is satisfied for all {\bf radial} functions $f$ supported in the ball $B_r(e)$ (and every $h\in \ell^2(\Gamma)$), then $(\Gamma,\cG)$ satisfies the property of {\bf radial rapid decay} with polynomial parameter ${\bf b}$. 
\item  If (\ref{RD-def-1}) is satisfied for the family of normalized spherical shell averages $\sigma\sigma_r$ with width $L\ge L_0>0$ (and every $h\in \ell^2(\Gamma)$), then $(\Gamma,\cG)$ satisfies the property of {\bf rapid decay of spherical shells},  with polynomial parameter ${\bf b}$. 
\end{enumerate} 
\end{definition}
Clearly, rapid decay implies radial rapid decay, with the same polynomial parameter. In turn, radial rapid decay implies rapid decay of spherical shells, with the same polynomial parameter.  As before, we suppress the width $L$ of the spherical shell from the notation, and note that under our assumptions below this will not affect the polynomial parameter $\bf b$, if $L$ is chosen sufficiently large.  

We note that the property of radial rapid decay is useful only when the length function takes a discrete set of values (for example in $a\NN$ for some $a > 0$). Indeed, consider such length functions as $D(\gamma p,p)$ or $\norm{\pi(\gamma)}$ associated with an invariant metric $D$ (on a symmetric space, for example) or a norm in a linear representation $\pi$. Then the sphere $S_r$ is empty except for a countable set of distances $r$, and $\abs{S_r}$ is a highly irregular function, so  that it is the spherical shells that must be considered in this case.

Let us now state the following property which will play a key role in our discussion, and which is motivated by the discussion in \cite{NT09}. 
\begin{definition}[rapid decay of spherical shell correlations] \label{Def-weak-RD}
Let $\Gamma$ be a countable group, and $\cG$ a length function. $\Gamma$ satisfies 
{\bf  rapid decay of spherical shell correlations} if there exist a non-negative constant ${\bf b}$ and a positive constant $L_0$ such that for any given  $L\ge L_0> 0$, and for arbitrary non-empty finite sets $A$ and $B$, and any $r \ge 1$  
\begin{equation}\label{weak-RD} 
\abs{\set{u\in A, v\in B\,;\, r\le d_\cG(u,v)< r+L}} \le C_{cor}  r^{\bf b} \sqrt{\abs{A}}\sqrt{\abs{B}}\sqrt{\abs{SS_r}}
\end{equation}
where $C_{cor}$ may depend on $L$, but is independent of $A$, $B$ and $r$. 
\end{definition}
We now note the following straightforward fact. 

\begin{prop}\label{RD-prop}
Rapid decay of spherical shells averages namely rough radial rapid decay, implies rapid decay of spherical shell correlations (with the same width $L$), with no change in the polynomial parameter. 
\end{prop}
\begin{proof}
Let us note that since $\cG(u^{-1}v)=d_\cG(u,v)$ we have 

\begin{equation}\label{dist-R-1}
\set{u\in A, v\in B\,;\, d_\cG(u,v)=r}=\bigcup_{w\in S_r}\set{u\in A, v\in B, u^{-1}v=w}
\end{equation}
and therefore 
\begin{equation}\label{dist-R-11}
\set{u\in A, v\in B\,;\, r\le d_\cG(u,v)< r+L}=\bigcup_{w\in SS_r}\set{u\in A, v\in B, u^{-1}v=w}.
\end{equation}
The various sets in this union are mutually disjoint, and using (\ref{conv-estim-1}) and the symmetry of $SS_r$, it follows that the number of elements in the set (\ref{dist-R-11}) is equal to   \begin{equation}\label{dist-R-2}\sum_{w\in SS_r} 1_{A^{-1}} \ast 1_B(w)=1_{A^{-1}}\ast 1_B \ast 1_{SS_r}(e).
\end{equation}

Furthermore, denote $f^\vee(\gamma)=f(\gamma^{-1})$.  Since for functions $f,h$ which are real-valued  $\inn{f^\vee,h}=\sum_{\gamma\in \Gamma} f(\gamma^{-1})h(\gamma)=f\ast h(e)$ we conclude 
\begin{equation}\label{dist-R-3}
1_{A^{-1}}\ast 1_B \ast 1_{SS_r}(e)= \inn{1_{A}, 1_B\ast 1_{SS_r}}=\abs{\set{u\in A, v\in B\,;\, r\le d_\cG(u,v)< r+L}}
\end{equation}

Assuming rapid decay of spherical shell averages, by (\ref{radial conv}) and  (\ref{RD-def-1}) 
\begin{equation}\label{sphere-conv}
\norm{\cS\cS_r (h)}_{\ell^2}\le C_{rd}\,r^{\bf b}\abs{SS_r}^{-1}\norm{1_{SS_r}}_{\ell^2}\norm{h}_{\ell^2}= C_{rd}\, r^{\bf b}\abs{SS_r}^{-1/2} \norm{h}_{\ell^2}, 
\end{equation} 

and using the Cauchy-Schwartz inequality we obtain the upper bound

\begin{equation}\label{dist-R-4}
\abs{ \inn{1_{A}, 1_B\ast 1_{SS_r}}}\le \norm{1_{A}}_{\ell^2}\abs{SS_r}\norm{\cS\cS_r(1_B)}_{\ell^2}\le  C_{rd}\,r^{\bf b} \sqrt{\abs{SS_r}}\sqrt{\abs{A}}\sqrt{\abs{B}}.
\end{equation}

We conclude that 
$$
\abs{\set{u\in A, v\in B\,;\, r\le d_\cG(u,v)<r+L}}\le C_{rd}\,r^{\bf b}\sqrt{\abs{A}}\sqrt{\abs{B}}\sqrt{\abs{SS_r}}\,.
$$

\end{proof}

{\bf Remarks.}
\begin{enumerate}
 \item In the case of the homogeneous tree of valency $2k$, the rapid decay of the spherical averages $\cS_r$ was shown by Cohen and Pytlik to be given by $\left(1+\frac{k-1}{k} r \right)(2k-1)^{-r/2}$, (see e.g. \cite{Ni17}\cite{Ni24} for further discussion). The rapid decay of spherical correlations is then also bounded by the same expression.
\item Let the length function $\cG$ assume only integer values. Then rapid decay of spherical shell averages implies radial rapid decay, but the polynomial parameter may be compromised. This is shown in \cite[Thm. 5.3] {Ni24} in the case the growth of spheres is almost exactly exponential, and the same proof applies if the growth of spheres is almost exactly polynomial-exponential. 
 \item In general, the polynomial parameter occurring in the rapid decay estimate of any group must be at least $1/2$, as shown by Nica \cite{Ni10}. 
 \item It was recently established by Boyer, Lobos and Pittet \cite{BLP22} that radial rapid decay for a word metric on a countable group does not imply rapid decay. This result will be discussed further in \S 5.2. 
\end{enumerate}

\subsection{Comparison bounds and the strong maximal inequality}

Let us consider the following straightforward comparison arguments. 
  
  \begin{prop}\label{integer-comparison}
   
  Let $(\Gamma, \cG)$ satisfy almost exact polynomial-exponential growth of spherical shells with parameters $(q,d)$ (and width $L\ge L_0> 0$). Then for $r \ge 1+r_0$ :
\begin{enumerate} 
\item The normalized averages $\sigma\sigma_r$ satisfy the inequality 
\begin{equation}\label{integers-1}
\sigma\sigma_r \le C( \sigma\sigma_{[r]}+ \sigma\sigma_{[r]+1}).
\end{equation} 
  \item  The balls also satisfy almost exact polynomial-exponential growth with parameters $(q, d)$, and the normalized averages satisfy 
  \begin{equation}\label{integers-2}\beta_r\le C(\beta_{[r]}+\beta_{[r]+1}).
  \end{equation}
  \item The Hardy-Littlewood maximal function is dominated by the supremum over balls with integer values : $\sup_{ r > 0}\beta_r\ast f\le C \sup_{j\in \NN} \beta_j\ast f$, for all $f \ge 0$. 
  \item If the spherical shells averages satisfy rapid decay with polynomial parameter $\bf b$, the ball averages satisfy rapid decay with polynomial parameter ${\bf b}+\frac12$. 
\end{enumerate}
  \end{prop}

\begin{proof}

(1) The inequality follows from the obvious inclusion of the supports of the two functions in (\ref{integers-1}), and the assumption of almost exact polynomial exponential growth of spherical shells.

(2) We have by the almost exact polynomial-exponential growth assumption with $q > 1$, for $r \ge 1+r_0+3L$  
$$r^dq^r\le C_1\abs{SS_{r}}\le C_2 \abs{SS_{r-L}} \le C_2\abs{B_r}\le C_3\abs{B_{[r]+1}}\le C_3\sum_{j=0}^{[r]+1} \abs{SS_j}.
$$
Now  for $k \ge 1$ we  have $C^{-1}_0 k^d q^{k}\le \sum_{j=1}^k j^d q ^j \le C_0 k^d q^{k}$. Therefore we obtain the bound    
$$
\le C_4 \sum_{j=0}^{[r]+1} (1+j)^d q^j\le C_4 ([r]+1)^{d}q^{[r]+1}\le C_5 r^d q^r\,.
$$
This shows that the balls have polynomial-exponential growth, and then inequality (\ref{integers-2}) is proved in the same manner as (\ref{integers-1}). 

(3) Since the balls have almost exact polynomial-exponential growth by (2), 
clearly $\beta_r \le C\beta_{[r]+1}$. Therefore 
$\sup_{r > 0}\beta_r \ast f \le C\sup_{j\in \NN}\beta_j\ast f$ for every $f \ge 0$.  
 
(4) Let us first restrict $r\ge 1+r_0+3L$, as well as $L$, to be integers. By definition $B_r= \cup_{j=0}^{r-L} SS_j$ and so 
\begin{equation} \label{compare}\beta_r=\frac{1_{B_r}}{\abs{B_r}}\le  \sum_{j=0}^{r-L}\mu_j\frac{1_{SS_j}}{\abs{SS_j}}=  \sum_{j=0}^{r-L}\mu_j\sigma\sigma_j\text{, where } \mu_j > 0\,\,,\,\, \sum_{j=0}^{r-L}\mu_j\le 1\,.
\end{equation} 
Here $\mu_j \le \frac{\abs{SS_j}}{\abs{B_r}}\le C_0q^{j-r}$ for $0\le j \le r-L$, since the  almost exact polynomial-exponential growth of balls and spherical shells is the same. Hence, since $\norm{\sigma\sigma_j}\le C_1 (j+1)^{\bf b} \abs{SS_j}^{-1/2}$, we conclude that 
$$\norm{\beta_r}_{\ell^2\to \ell^2}\le \sum_{j=0}^{r-L}\mu_j \norm{\sigma\sigma_j}_{\ell^2\to \ell^2}\le \sum_{j=0}^{r-L}C_0 q^{j-r} C_1 (j+1)^{\bf b} \left(C_2(j+1)^d q^{j}\right)^{-1/2}$$
$$
\le C_3 q^{-r} \sum_{j=0}^{r}  (j+1)^{\bf b}q^{j} \left((j+1)^d q^{j}\right)^{-1/2}$$
$$\le 
C_3q^{-r} \cdot\left(\sum_{j=0}^{r} (j+1)^{2{\bf b}}\right)^{1/2} \left(\sum_{j=0}^{r}q^{2j}  \left((j+1)^d q^{j}\right)^{-1}\right)^{1/2}
$$
using Cauchy-Schwarz.  The two factors in the foregoing inequality are estimated by $\sum_{j=0}^r (j+1)^{2{\bf b}}\le C_6 r^{2{\bf b}+1}$, and 
$\sum_{j=0}^r (j+1)^{-d} q^{j}\le C_7 r^{-d}q^{r}$, where $C_6$ and $ C_7$ depends on $q$ and $d$ but not on $r$. For the second estimate it suffices to note that (setting $i=r-j$) for $r \ge 2$ 
\begin{equation}\label{summation} \sum_{j=0}^r \left(\frac{r}{j+1}\right)^{d} q^{j-r}=\sum_{i=0}^{r}\left(1-\frac{i-1}{r}\right)^{-d} \left(\frac{1}{q}\right)^i\le C \sum_{i=0}^r (i+1)^d \left(\frac{1}{q}\right)^i < \infty.
\end{equation}
%

The bound we are interested in is 
$$
\norm{\beta_r}_{\ell^2\to \ell^2}\le C_4 r^{{\bf b}+\frac12} (r^d q^r)^{-1/2}\le C_5 r^{{\bf b}+\frac12}\abs{B_r}^{-1/2} \,, $$
and it follows easily using the foregoing estimates of the two factors.. 

Finally by (\ref{integers-2}) we conclude that for any $r\ge 1+r_0$ and $f \ge 0$
$$\norm{\beta_r\ast f}_{\ell^2}\le C_6( \norm{  \beta_{[r]}\ast f}+ \norm{  \beta_{[r]+1}\ast f})\le C_8 r^{{\bf b}+\frac12} (r^d q^r)^{-1/2}\norm{f}_{\ell^2} 
$$
and since $\norm{\beta_r\ast f}_{\ell^2}\le \norm{\beta_r\ast \abs{f}}_{\ell^2}$ we obtain the stated bound for any ball averaging operator $\beta_r$, $r \ge 1$. 

\end{proof} 
It is now straightforward to deduce the following result regarding the strong maximal inequality in $\ell^p(\Gamma)$, $1< p \le  \infty$. 
\begin{prop}
Let $(\Gamma, \cG)$ satisfy almost exact polynomial-exponential growth for spherical shells as well as rapid decay of the spherical shell averages.
Then the  Hardy-Littlewood maximal operators $\cM$ defined in (\ref{def-HL-max-op}) satisfies the strong maximal inequality in $\ell^p$, $1 < p \le \infty$, namely $\norm{\cM f}_{\ell^p}\le C_p\norm{f}_{\ell^p}$. 
\end{prop}
\begin{proof}
By (\ref{integers-2}) to bound the Hardy-Littlewood maximal operator for general radius $r$ it suffices to bound $\cM_{I}=\sup_{j\in \NN }\beta_j$.
 The averages $\beta_j$ satisfy the same two properties as $\sigma\sigma_r$ by Proposition \ref{integer-comparison}.  
Consequently $\norm{\beta_j}\le C_0 e^{-\alpha j}$, where we can take for example  $\alpha=\frac14\log q > 0$. 
Then for any $ f\in \ell^2$, 
$$ \sum_{j=0}^\infty e^{j\alpha/2}\norm{\beta_j\ast \abs{f}}_{\ell^2}\le C_1 \norm{f}_{\ell^2}\,.$$
This amounts to an effective exponential strong maximal inequality in $\ell^2$, 
which immediately implies the strong type $(2,2)$-maximal inequality for $\beta_j$, $j\in \NN$, and hence for the family $\beta_r$. The case of the strong type $(p,p)$-maximal inequality for $1 < p < \infty$ is obtained similarly. 
Of course, the same considerations apply also to the family $\sigma\sigma_r$, and in fact to any sequence of averages in any group algebra whose convolution norms decay exponentially. 

The strong maximal inequality for $p=\infty$ is obvious by definition, but an exponential decay estimate does not hold, as $\norm{\sigma\sigma_r}_{\ell^\infty(\Gamma)\to \ell^\infty(\Gamma)}=1$ for all $r$. 
\end{proof}

\section{Distributional inequality for spherical shell averages}  
We now turn to state and prove a distributional inequality for each of the operators $\cS \cS_r$ supported on the spherical shells $SS_r$. Choosing $L_0$ suitably we can and will assume that $\abs{SS_r}\ge 1$ for all $r\ge 0$, and consider any radius $r \ge 1$. The statement and proof of the following proposition are motivated by a remarkable argument developed in \cite{NT09} for the rooted or regular tree. The argument in question uses radial considerations only, and the following result constitutes a direct generalization of it. 

\begin{prop}\label{distribution-estm}
Assume $(\Gamma, \cG)$ satisfies rapid decay of spherical shell correlations (\ref{weak-RD}) with polynomial parameter ${\bf b}$. Then there exists $C_{dis}> 0$ such that 
for any non-negative $f\in \ell^1(\Gamma)$ the distribution of the spherical shells averages $\cS\cS_r$, $r\ge 1$ satisfies, $\forall \eta > 0$:
\begin{equation}\label{max-ineq-proof-1}
\abs{\set{w\,;\, \cS\cS_r (f)(w)\ge \eta } }\le 
C_{dis} r^{2{\bf b}}
\sum_{1\le 2^n \le 2\abs{SS_r}} \sqrt{\frac{2^n}{\abs{SS_r}}} \cdot 2^n \abs{\set{w\,;\, f(w)\ge 2^{n-1}\eta}}
\end{equation}
\end{prop}
\begin{proof}
We assume $f$ to be non-negative, and replacing it by $f/\eta$, we can assume $\eta = 1$. 
Consider the level sets of $f$ given by  
$$A_n=\set{w\,;\, 2^{n-1}\le f(w)\ < 2^n}\,\,,\,\,n \ge 0$$
 and note that 
\begin{equation}\label{max-ineq-proof-2} 
f\le \frac12 +\sum_{1\le 2^n\le \abs{SS_r}}2^n 1_{A_n}+f1_{\set{f \ge \frac12\abs{SS_r}}}
\end{equation} 
Averaging $f$ on a spherical shell of radius $r$ we deduce:
\begin{equation}\label{max-ineq-proof-3}
\cS\cS_r (f)\le \frac12+\sum_{1\le 2^n\le \abs{SS_r}}2^n \cS\cS_r (1_{A_n})+
\cS\cS_r \left(f1_{\set{f \ge \frac12  \abs{SS_r}}}\right)
\end{equation}  
To estimate this expression, we start with the third summand. Assuming that $\cS\cS_r \left(f1_{\set{f \ge \frac12 \abs{SS_r}}}\right)(x)\neq 0$ at some point $x$, it follows that  the distance of $x$ from the set $\set{y\,;\, f(y) \ge \frac12 \abs{SS_r}}$ is smaller than $r+L$. 
Hence the support of the function $\cS\cS_r \left(f1_{\set{f \ge \frac12  \abs{SS_r}}}\right)$ is contained in the union 
of the spherical shells $SS_r$ centered at the points of the set  $\set{y\,;\, f(y) \ge \frac12 \abs{SS_r}}$. Therefore
\begin{equation}\label{max-ineq-proof-4} 
\abs{\set{x\,;\, \cS\cS_r \left(f1_{\set{f \ge \frac12  \abs{SS_r}}}\right)(x)\neq 0}} \le \abs{SS_r}\cdot \abs{\set{f \ge \frac12  \abs{SS_r}}}\,.
\end{equation} 
We now claim  that the size of the set $\set{\cS\cS_r f\ge 1 } $ is bounded by 
\begin{equation}\label{max-ineq-proof-5}  
 \abs{\set{x\,;\,  \sum_{1\le 2^n\le \abs{SS_r}}2^n \cS\cS_r (1_{A_n})(x)\ge \frac12}}+ \abs{SS_r}\cdot \abs{\set{ f \ge \frac12 \abs{SS_r}}}\,.
\end{equation}
Indeed, $\cS\cS_r f(x)\ge 1 $ implies using (\ref{max-ineq-proof-3}) that
$$\frac12 \le \sum_{1\le 2^n\le \abs{SS_r}}2^n \cS\cS_r (1_{A_n})(x)+
\cS\cS_r \left(1_{\set{f \ge \frac12  \abs{SS_r}}}\right)(x)$$
and this implies that if the second summand in this estimate vanishes, necessarily $\sum_{1\le 2^n\le \abs{SS_r}}2^n \cS\cS_r (1_{A_n})(x)\ge \frac12 $. 
Therefore the size of the set $\set{\cS\cS_r f\ge 1 } $ is bounded by the sum of the size of the support of $\cS\cS_r \left(1_{\set{f \ge \frac12  \abs{SS_r}}}\right)$ (which is bounded 
by $\abs{SS_r}\cdot \abs{\set{ f \ge \frac12 \abs{SS_r}}}$ as we just saw) together with the size of the set $$\set{x\,;\,  \sum_{1\le 2^n\le \abs{SS_r}}2^n \cS\cS_r (1_{A_n})(x)\ge \frac12}$$ as claimed in  (\ref{max-ineq-proof-5}).

Consider now the first summand of (\ref{max-ineq-proof-5}).
We claim that if 
$$ \sum_{1\le 2^n\le \abs{SS_r}}2^n \cS\cS_r (1_{A_n})(x)\ge \frac12$$ 
then for at least one $n$ with $1\le 2^n \le \abs{SS_r}$ 
$$\cS\cS_r (1_{A_n})(x)\ge \frac{1}{2^{n+4}} \left(\frac{2^n}{\abs{SS_r}}\right)^{1/4}.$$
Indeed, otherwise when summing the geometric series $2^{n/4}$ on the interval
$1\le 2^n\le \abs{SS_r}$, we obtain the contradiction 
$$\frac12 \le  \sum_{1\le 2^n\le \abs{SS_r}}2^n \cS\cS_r (1_{A_n})(x)\le 
\frac{1}{16}\sum_{ 1\le 2^n\le \abs{SS_r}} \left(\frac{2^n}{\abs{SS_r}}\right)^{1/4}
$$
$$
=\frac{1}{16}\frac{1}{\abs{SS_r}^{1/4}}\sum_{n=0}^{\log_2 \abs{SS_r}} \left(2^{1/4}\right)^n\le \frac{2^{1/4}\abs{SS_r}^{1/4}-1}{16\cdot \abs{SS_r}^{1/4}(2^{1/4}-1)}
$$
$$ 
\le 
\frac{2^{1/4}\abs{SS_r}^{1/4}}{16\cdot \abs{SS_r}^{1/4}(2^{1/4}-1)}
<\frac12.
$$

Define now 
\begin{equation}\label{max-ineq-proof-6} 
B_n=\set{x\,;\, \cS\cS_r (1_{A_n})(x)\ge \frac{1}{2^{n+4}} \left(\frac{2^n}{\abs{SS_r}}\right)^{1/4}}.
\end{equation} 
Then by (\ref{max-ineq-proof-5}) and the argument following it, 
$\set{\cS\cS_r f\ge 1 } $ is bounded by 
\begin{equation}\label{max-ineq-proof-7}   
\set{\cS\cS_r f\ge 1 }\le \sum_{1\le 2^n \le \abs{SS_r}}\abs{B_n} +\abs{SS_r}\cdot \abs{\set{ f \ge \frac12 \abs{SS_r}}}
\end{equation}
To further estimate this bound, note that using (\ref{dist-R-3})  and the definition of  $B_n$ 
$$\frac{1}{\abs{SS_r}}\abs{\set{(u,v)\in B_n\times A_n\,;\, r\le d_\cG(u,v)< r+L}}$$
$$=\inn{1_{B_n}, \cS\cS_r(1_{A_n})}=\sum_{u\in B_n}\cS\cS_r(1_{A_n})(u)
\ge \frac{\abs{B_n}}{2^{n+4}}\left(\frac{2^n}{\abs{SS_r}}\right)^{1/4}$$

Now using our assumption that the bound (\ref{weak-RD}) holds, namely

$$\abs{\set{(u,v)\in B_n\times A_n\,;\, r\le d_\cG(u,v)< r+L}}\le 
 C_{cor}r^{{\bf b}} \sqrt{\abs{A_n}}\sqrt{\abs{B_n}}\sqrt{\abs{SS_r}}
$$
we conclude that 
$$
\frac{\abs{B_n}}{2^{n+4}}\left(\frac{2^n}{\abs{SS_r}}\right)^{1/4}\le 
C_{cor}r^{{\bf b}} \sqrt{\frac{\abs{A_n}\abs{B_n}}{\abs{SS_r}}}
$$
so that squaring the last inequality 
$$
\abs{B_n}\le \frac{(C_{cor}r^{{\bf b}})^2\abs{A_n}}{\sqrt{\abs{SS_r}}}2^{2n+8-\frac12 n}=2^8C_{cor}^2r^{2{\bf b}}\sqrt{\frac{2^n}{\abs{SS_r}}} \cdot 2^n \abs{A_n}\,.
$$
Using the bound just found for $B_n$ together with (\ref{max-ineq-proof-7}) we conclude that the size of the set $\set{\cS\cS_r f\ge 1 }$ is bounded by 
\begin{equation}\label{final sum}
C_{0} r^{2{\bf b}}\sum_{1\le 2^n \le \abs{SS_r}} \sqrt{\frac{2^n}{\abs{SS_r}}} \cdot 2^n \abs{A_n}
 +\abs{SS_r}\cdot \abs{\set{ f \ge \frac12 \abs{SS_r}}}
\end{equation} 
with $C_{0}$ chosen to satisfy $C_{0} \ge 1$.  
Now let $n_l\ge 1$ be the last poitive integer such that $2^{n_l} \le 2\abs{SS_r} $. 
Then  $\abs{SS_r} \le  2^{n_l} $ so that  

$$
2^{n_l-1}\le \abs{SS_r}\le 
2^{n_l}\text{,  and  }  \frac{1}{\sqrt{2}} \le  \sqrt{\frac{2^{n_l-1}}{\abs{SS_r}}}\le 1  \,. $$
The second summand in (\ref{final sum}) can then be estimated by
$$
   2 \sqrt{ \frac{2^{n_l-1}}{\abs{SS_r}}}\cdot 2^{n_l-1} \abs{\set{ f \ge 2^{n_l-2}}},
$$
 since $\set{f \ge \frac12 \abs{SS_r}}\subset \set{f \ge 2^{n_l-2}}$.

Finally, we note that the sum (\ref{final sum}) can be estimated by a single formula, as stated in Proposition \ref{distribution-estm}. Recall $A_n=\set{w\,;\, 2^{n-1}\le f(w)\ < 2^n}\,\,,\,\,n \ge 0$ and multiply the second summand of (\ref{final sum}) by  $2C_0  r^{2{\bf b}}$. Now  let $2^n$ range over the larger interval $[1, 2\abs{SS_r}]$.  We then conclude that  for $r\ge 1$, setting $C_{dis}=2C_0$ 
$$
  \set{\cS\cS_r f\ge 1 }\le C_{dis} r^{2{\bf b}}
\sum_{1\le 2^n \le 2\abs{SS_r}} \sqrt{\frac{2^n}{\abs{SS_r}}} \cdot 2^n \abs{\set{w\,;\, f(w)\ge 2^{n-1}}}$$
 
This completes the proof of Proposition \ref{distribution-estm}. 
\end{proof}

\section{Proof of the maximal inequality for ball averages} 
\subsection{Weak-type maximal inequality in Orlicz space}\label{Sec: orlicz}

 Let $(X, \cB,\mu)$ be a $\sigma$-finite measure space with measure $\mu$, which we assume is standard, for example isomorphic to $\NN$ with the counting measure, or  $[0,1]$ or $\RR$ with Lebesgue measure. Following \cite[p.275]{Fa72} and \cite[\S 3]{Fa84} define the function space $\cL\left(\log \cL\right)^{\bf c}$ (where ${\bf c} > 0$) consisting of measurable functions $f$ on $X$  such that  for all $ \eta > 0$
\begin{equation}\label{Orlicz} 
\int_{\set{\abs{f} >  \eta}} \frac{\abs{f(x)}}{\eta}\left(\log \frac{\abs{f(x)}}{\eta}\right)^{\bf c} d\mu(x) < \infty.
\end{equation} 
Note that since ${\bf c} >  0$ 
$$\int_{\set{\abs{f}>  \eta}} \frac{\abs{f(x)}}{\eta}\left(\log \frac{\abs{f(x)}}{\eta}\right)^{\bf c} d\mu(x) =\int_X \frac{\abs{f(x)}}{\eta}\left(\log^+ \frac{\abs{f(x)}}{\eta}\right)^{\bf c} d\mu(x)$$
$$=\int_{\set{\abs{f}\ge   \eta}} \frac{\abs{f(x)}}{\eta}\left(\log \frac{\abs{f(x)}}{\eta}\right)^{\bf c} d\mu(x) 
$$
Furthermore $f$ is in the usual Lebesgue space $ L\left(\log^+ L\right)^{\bf c}$ if and only if the integral (\ref{Orlicz})
is finite for at least one $\eta > 0$. It follows that $\cL\left(\log \cL\right)^{\bf c}\subset L\left(\log^+ L\right)^{\bf c}$. 

The function space has the following properties. 
\begin{enumerate}
\item $\cL\left(\log \cL\right)^{\bf c}$ is a vector subspace of $L\left(\log^+ L\right)^{\bf c}\subset L^1+L^\infty$. 
\item
$\cL\left(\log \cL\right)^{\bf c}= L\left(\log^+ L\right)^{\bf c}$ if and only if $\mu(X)< \infty$. 
\item $\bigcup_{1< p < \infty} L^p \subsetneq \cL\left(\log \cL\right)^{\bf c}$. 
\item When ${\bf c} \ge 1$, $L\left(\log^+ L\right)^{\bf c}$ is a Banach space under the norm 
$$\norm{f}=\inf\set{\eta > 0\,;\, \int_{\set{\abs{f}>  \eta}} \frac{\abs{f(x)}}{\eta}\left(\log \frac{\abs{f(x)}}{\eta}\right)^{\bf c} d\mu(x)\le 1}$$
and $\cL\left(\log \cL\right)^{\bf c}$ is a closed subspace. 
For $f\in \cL\left(\log \cL\right)^{\bf c}$ with ${\bf c}\ge 1$ we will denote this norm 
by $\norm{f}_{L\left(\log^+ L\right)^{\bf c}}$. 
\end{enumerate} 

\subsection{Weak-type maximal inequalities} 
Recall that the operator of averaging a function on a ball of radius $r$ is denoted by 
$\cB_r f=f\ast \beta_r$.  The maximal operator we consider is (the global part of) the Hardy-Littlewood maximal operator, given by 
\begin{equation}\label{HL-op}
\cM f(x)=\sup_{r >  r_0} \cB_r \abs{f}(x).
\end{equation}

\begin{definition}
The Hardy-Littlewood operator $\cM$ satisfies the weak-type maximal inequality in $\cL\left(\log \cL \right)^{\bf c}(X) $ if there exists a positive constant $K_{\cM}$ (independent of $f$), such that for every $\eta > 0$ 
$$\abs{\set{x\,;\, \cM f(x)\ge \eta}}\le  K_{\cM}\sum_{\abs{f} > \eta}  \frac{\abs{f(x)}}{\eta} \left(\log \frac{\abs{f(x)}}{\eta}\right)^{\bf c}
$$
\end{definition}

Our standing assumption is that $SS_r$, $r \ge 1$  has almost exact polynomial-exponential growth, with $q > 1$, see (\ref{exact growth}). Then by Proposition \ref{integer-comparison}  
the maximal function associated with the spherical shell averages $\cS\cS_r$ with integer radius  $r\in \NN_+$ dominates the maximal function associated with $\cS\cS_r$, $r\ge 1$. Similarly, the maximal function associated with ball averages $\cB_r$ with integer radius $r\in \NN_+$ dominates the maximal function associated with $\cB_r, r\ge 1$. Furthermore, using (\ref{compare}) clearly the operator $\cB_r$ for $r\in \NN_+$ is dominated by a convex combination of the operators $\cS\cS_s$, $s\in \NN_+$ with $1\le s \le r$, together with $\cB_L=\cS\cS_0$. Since $\cB_L: \ell^1(\Gamma)\to \ell^1(\Gamma)$  is a fixed averaging operator with norm $1$, it follows from standard considerations that it suffices to show that $\cM$ is dominated by the maximal operator associated with spherical shell averages $\cS\cS_r, r\in \NN_+$, namely by 
\begin{equation}\label{sphere-max-op}
\cA f(x)=\sup_{r\in \NN_+} \cS\cS_r \abs{f}(x).
\end{equation}

We now turn to state and prove our main result, the weak-type maximal inequality for ball averages. The proof is based on the distributional inequality stated in Proposition \ref{distribution-estm}, and is motivated by another important argument developed in \cite{NT09} for the rooted or regular tree. 

\begin{thm}[Main Theorem : Rapid Decay of Spherical Shell Correlations]\label{main}
Let $(\Gamma,\cG)$ satisfy the property of rapid decay of spherical shell correlations (\ref{weak-RD}) with polynomial parameter ${\bf b}$, and satisfy almost exact exponential-polynomial growth of spherical shells (with any fixed width $L\ge L_0>0$). 
Then the Hardy-Littlewood operator $\cM$ for ball averages satisfies the weak-type  $\cL\left(\log \cL\right)^{\bf c}(\Gamma)$ maximal inequality, where ${\bf c}=2{\bf b}$.
\end{thm} 

\begin{proof}
We can and will assume that $f \ge 0$. 
As noted in the previous paragraph, it suffices to bound the operator $\cA$,  
which clearly satisfies :
$$\abs{\set{\cA  f\ge  \eta}}\le C\sum_{r\in \NN_+} \abs{\set{\cS\cS_r f \ge \eta}}\,.$$
We now turn to bound the right hand side. 

Using Proposition \ref{distribution-estm}, and the two-sided bounds for the growth of $\abs{SS_r}$ stated in (\ref{exact growth}), we have :
$$\sum_{r\ge 1} \abs{\set{\cS \cS_r f \ge \eta}}\le
$$
$$C_{dis}\sum_{r\ge 1}  r^{2{\bf b}} \sum_{1\le  2^n  \le 2\abs{SS_r}}
\sqrt{\frac{2^n}{\abs{SS_r}}}\cdot 2^n \abs{\set{w\,;\, f(w)\ge 2^{n-1} \eta}}
$$
$$
= C_{dis}\sum_{n \in \NN} \left(\sum_{r\ge 1\,;\, 2\abs{SS_r}\ge 2^{n}} \frac{ r^{2{\bf b}}}{\sqrt{\abs{SS_r}}}\right) 2^{3n/2} \abs{\set{w\,;\, f(w)\ge 2^{n-1} \eta}}$$

$$\le C_{dis} \sum_{n \in \NN} \left(\sum_{2C_{gr} r^{d}q^r\ge 2^{n}} \frac{C^{1/2}_{gr} r^{2{\bf b}}}{\sqrt{r^{d}q^r}}\right) 2^{3n/2}\abs{\set{w\,;\, f(w)\ge 2^{n-1} \eta}} \le $$

\begin{equation}\label{large ineq}
 C_{dis}C_{gr}^{1/2} 
 \sum_{n \in \NN}
  \left(\sum_{(2C_{gr})^{1/2}  r^{d/2}q^{r/2}\ge 2^{n/2} }
\frac{r^{2{\bf b} -\frac{d}{2}}}{q^{r/2}}\right)2^{3n/2}\abs{\set{w\,;\, f(w)\ge 2^{n-1} \eta}}. 
\end{equation} 
To estimate (\ref{large ineq}), let us more generally fix $s > 1$, $ a \ge 0$, $b\in \RR$, $ c \ge 1$, and $\Upsilon \ge 1$.  Let $k_0$  be the least positive integer such that $ck_0^{a}s^{k_0}\ge  \Upsilon $. 
 
Then $k_0\ge 1$ and $ck^{a} s^{k}\ge  \Upsilon$ for all $k \ge k_0$. Setting $k-k_0=m$, $\frac{k}{k_0}=1+\frac{m}{k_0}$ we have
$$\sum_{k\,;\,ck^{a} s^k \ge  \Upsilon} \frac{k^b}{s^k} = \sum_{k \ge  k_0}\frac{k^b}{s^k} 
\le \frac{k_0^{b}}{s^{k_0}}\sum_{k = k_0}^\infty \frac{\left(\frac{k}{k_0}\right)^b}{s^{k-k_0} }
$$
$$
= \frac{k_0^{b}}{s^{k_0}} \sum_{m=0}^\infty \frac{\left(1+\frac{m}{k_0}\right)^b}{s^m }
\le  \frac{k_0^{b}}{s^{k_0}} \sum_{m=0}^\infty \frac{\left(1+m\right)^b}{s^m }
\le  \lambda_{b,s} \frac{k_0^{b+a}}{ck_0^{a}s^{k_0}} \le  \frac{\lambda_{b,s} }{\Upsilon}k_0^{b+a}
$$ 
The last inequality in the previous line follows since $ck_0^{a} s^{k_0}\ge  \Upsilon$. 

We now claim that under the additional condition $a+b \ge 0$ the previous expression is bounded by

\begin{equation}\label{log est} 
  C_1\frac{\left(\log (\Upsilon+C_0) \right)^{b+a}}{\Upsilon} \,.
\end{equation} 
with $C_0> 0$ , $C_1 > 0$, depending on $c,s,a+b$. 

To prove this bound,  we estimate $k_0^{b+a}$, as follows. If  $k_0\le   3$ (say) 
 then clearly 
$1\le \Upsilon\le cs^{k_0} k_0^{b+a}\le cs^3 3^{a+b}$ so that $k_0^{b+a}\le c^{-1} s^23^{a+b}$, which is constant depending on $c,s,a+b$. So we can certainly write 
$k_0^{a+b} \le C_1\left( \log (\Upsilon+C_0)\right)^{a+b}$ for some $C_0> 0$ , $C_1 > 0$ (depending on $c,s,a+b$). 

Otherwise for $k_0 \ge 4$ note that by definition of $k_0$ necessarily $ c3^{a+b}  < c(k_0-1)^{a} s^{k_0-1}< \Upsilon $. Therefore can conclude that 
$$\frac12 k_0\log s\le \log c +(k_0-1)\log s +a\log (k_0-1)< \log \Upsilon$$
and consequently   
$$k_0^{b+a}\le \left(\frac{2}{\log s} \right)^{b+a}\log^{a+b}\Upsilon\,.$$

Now going back to (\ref{large ineq}), set $s=q^{1/2}> 1$, $a=\frac12 d\ge 0$,  $b=2{\bf b} -\frac12 d$ (so that $a+b=2{\bf b}\ge 0$), $c=(2C_{gr})^{1/2}$  and $\Upsilon= 2^{n/2}\ge 1 $ we conclude using (\ref{log est}) that 
$$ \sum_{(2C_{gr})^{1/2}r^{d/2}q^{r/2}\ge 
2^{n/2} } \frac{r^{2{\bf b}-\frac{d}{2}}}{q^{r/2}}\le C_2\frac{ 1}
{2^{n/2}}\left(\log(2^n+C_3)\right)^{a+b}$$
$$\le C_4 \cdot 2^{-n/2}(n+1)^{a+b}.
$$
with $C_4$ depending on $(d,q, {\bf b}, C_{gr})$ but independent of $r$ and $n$. 

Noting that  $a+b=2{\bf b}$, using the foregoing estimate we conclude that the maximal inequality holds in $\cL\left(\log \cL\right)^{2{\bf b}}(\Gamma)$, as follows. 
 
$$\sum_{r\ge 1} \abs{\set{\cS\cS_r f\ge \eta}}=\sum_{r\ge 1} \abs{\set{\cS\cS_r \left(\frac{f}{\eta} \right)\ge 1}}
$$
$$
\le 
C_5\sum_{n \in \NN}2^{-n/2} 2^{3n/2}(n+1)^{2{\bf b}} \abs{\set{x\,;\, \frac{f(x)}{\eta} \ge 2^{n-1}}}
$$
$$
=
 C_5 \sum_{x \in \Gamma}\sum_{n=0}^\infty 2^n (n+1)^{2{\bf b}}\cdot 1_{\set{x\,;\, \frac{f(x)}{\eta} \ge 2^{n-1}}}$$
 $$\le C_6\sum_{x \in \Gamma}\left(2+\left(\sum_{n=1}^{\lfloor\log^+_2\left(\frac{f(x)}{\eta}\right)\rfloor}   2^n (n+1)^{2{\bf b}}\right)\right)\,
$$
where the upper limit of the summation satisfies $\lfloor\log^+_2\left(\frac{f(x)}{\eta}\right)\rfloor\ge 1$ or equivalently $\frac{f(x)}{\eta}\ge 2$. 

We now use the fact that for $ N=\lfloor\log^+_2\left(\frac{f(x)}{\eta}\right)\rfloor \ge 1$ 
we have (using an argument similar to (\ref{summation})) : 
$$2+ \sum_{n=0}^N 2^n (n+1)^{2{\bf b}}\le C_7 2^N  N^{2{\bf b}}\,,$$ 
together with the fact that $2^N \le\frac{f(x)}{\eta}$ and the fact that ${\bf c}:=2{\bf b} \ge 0$, 

and conclude that the following bound holds : 

$$
  \sum_{r\ge 1} \abs{\set{\cS\cS_r f\ge \eta}}\le K\sum_{x\in \Gamma} \frac{f(x)}{\eta}\left(\log^+\frac{f(x)}{\eta}\right)^{2{\bf b}} 
$$

and this completes the proof of Theorem \ref{main}. 
\end{proof}

\section{Weak-type maximal inequality for lattice subgroups} 

\subsection{Lattice points, almost exact growth and radial RD} 
Let $G$ be a non-compact lcsc group, and let $\Gamma$ be a discrete countable lattice subgroup. It follows that $G$ is a unimodular group.

{\it Standing assumption: admissible length functions}. The length function $\tilde{\cG}$ on $G$ will be assumed throughout this section to be non-negative, proper and locally bounded, symmetric and subadditive, Borel measurable but not necessarily continuous, and satisfy $\tilde{\cG}(e)=0$. 
 We assume that for some positive $r_0$ the measure of a ball of radius $r \ge  r_0> 0 $ is positive. Such length functions will be called admissible. 
 
The restriction of $\tilde{G}$ to $\Gamma$ will be  denote by $\cG$.

\begin{definition}\label{lattice point count} 
The lattice $\Gamma$ has 
\begin{enumerate}
\item an effective solution to the lattice point counting problem in the balls $\tilde{B}_r$ if there exist a positive constant $ \kappa$ so that 
\begin{equation}\label{effective count}
\frac{\abs{\tilde{B}_r\cap \Gamma}}{m_G(\tilde{B}_r)}=1+O\left(m_G(\tilde{B}_r)^{-\kappa}\right)\,\,,\,\, \forall r \ge 1\,,
\end{equation}
\item an exact solution to the lattice point counting problem in the balls $\tilde{B}_r$ if
\begin{equation}\label{exact count}
\lim_{r \to \infty}  \frac{\abs{\tilde{B}_r\cap \Gamma}}{ m_G(\tilde{B}_r)} =1\,,
\end{equation}
\item  an almost exact solution to the lattice point counting problem in the balls $\tilde{B}_r$ if there exists a positive constant $C$ so that 
\begin{equation}\label{almost exact count}  C^{-1} m_G(\tilde{B}_r)\le \abs{\tilde{B}_r\cap \Gamma}\le Cm_G(\tilde{B}_r)\,\,,\,\, \forall r \ge r_0\,.
\end{equation}
\end{enumerate}
\end{definition} 
As in the discrete case, the balls $\tilde{B}_r=\set{g\in G\,;\,\tilde{ \cG}(g)< r}$, $r \ge  r_0> 0$  are said to have almost exact polynomial-exponential volume growth under (left and right) Haar measure $m_G$, if $C^{-1} r^d q^r \le m_G(\tilde{B}_r) \le C r^d q^r\,$ (for some fixed $C$ and for all $r \ge  r_0$). The condition of almost exact polynomial-exponential growth of spherical shell $\tilde{S}\tilde{S}_r$, $r\ge 1$ is the same, where the width satisfies $\tilde{L}\ge \tilde{L_0} > 0$. 

Let us note the following straightforward observations. 
\begin{lem}\label{restricted-gauge-volume}
\begin{enumerate}
Assume that the balls $\tilde{B}_r$ in $G$ have almost-exact polynomial-exponential growth w.r.t. the length function $\tilde{\cG}$.
\item Assume $\Gamma$ has almost exact solution to the lattice point counting problem in the balls $\tilde{B}_r$. Then the balls $B_r$ in $\Gamma$ w.r.t. the length function $\cG$ have almost exact polynomial-exponential growth with the same parameters $(d,q)$. 
\item Assume the balls $B_r\subset \Gamma$ w.r.t. the length function $\cG$ have almost exact polynomial-exponential volume growth. Then for $L\ge L_0$ ($L_0$ determined below), the spherical shells of width $L$ w.r.t. $\cG$ also have almost exact polynomial-exponential volume growth with the same parameters $(d,q)$.
 
\end{enumerate}
\end{lem} 
\begin{proof}
Part 1) is clear. For part 2), note that since 
$SS_r =B_{r+L} \setminus B_r$ we can estimate 
$$\frac{1}{C} (r+L)^d q^{r+L}-Cr^d q^r  \le \abs{SS_r}\le C(r+L)^d q^{r+L}-\frac{1}{C}r^d q^r$$
or equivalently 
$$\frac{1}{C} (1+\frac{L}{r})^dq^L-C\le \frac {\abs{SS_r}}{r^dq^r} \le C (1+\frac{L}{r})^dq^L-\frac{1}{C}$$
so that for $r \ge L$, and since $C \ge 1$
$$\frac{1}{C}q^L-C\le  \frac {\abs{SS_r}}{r^dq^r}\le C2^dq^L+1\,.$$ 
For every choice of $L \ge L_0=\frac{2 \ln C +\ln 2 }{\ln q}$ we have $\frac{1}{C}q^L-C\ge C > 0$ and the claim follows. 
\end{proof}

We now turn to discuss the radial rapid decay property of spherical shells in $\Gamma$, defined using a length function $\cG$ arising as the restriction of a length function $\tilde{\cG}$ defined on $G$.  
We will not require that  
$G$ itself satisfies the rapid decay property in $L^2(G)$, only that the spherical shells $\tilde{S}\tilde{S}_r$ have that property. The most obvious source of such groups are those that satisfy the property of rapid decay. Connected  Lie groups satisfying rapid decay have been characterized by Chatterji, Pittet and Saloff-Coste \cite[Thm. 0.1]{CPSC} who showed that the Lie algebra of $G$ is a direct sum of a semisimple Lie algebra and a Lie algebra of type $R$ (recall that the connected Lie group associated with a type $R$ Lie algebra has polynomial volume growth).  The growth of $G$ is therefore exponential if and only if the semisimple component of its Lie algebra is non-compact.  
Rapid decay holds then in groups somewhat more general than connected reductive non-compact Lie groups with compact center.  It is therefore a condition strictly weaker than the validity of the Kunze-Stein phenomenon, since the latter has been shown to hold for reductive groups with compact center and only for them, see \cite[Cor. 7.2]{Co78}. For discrete groups, the property of radial rapid decay is strictly weaker than rapid decay, as we will explain  in more detail below. We note that a general characterization of the connected Lie groups, or the algebraic groups over local fields, that satisfy radial rapid decay seems to be an open problem.

Our purpose is to exploit the properties of a family 
$\cU_t$ of sets in the lcsc group $G$ in order to deduce the rapid decay 
of the averages uniformly distributed on characteristic functions of the sets $\cU_t\cap \Gamma$ in $\ell^2(\Gamma)$. The properties we will utilize are almost-exact polynomial-exponential growth of $\cU_t$, rapid decay of the convolution norms of $\chi_{\cU_t}/m_G(\cU_t)$ on $L^2(G)$, almost exact solution of the lattice point counting problem for lattice points in $\cU_t$, and coarse stability of the sets $\cU_t$ under left and right perturbations. This last property which was introduced and discussed in \cite[\S 3.3]{GN10} is crucial in utilizing coarse comparison arguments between $\cU_t$ and $\cU_t\cap \Gamma$, and is defined as follows. 

\begin{definition}\label{coarse-def}{\bf Coarse admissibility.} 
{\rm
Let $G$ be an lcsc unimodular group with Haar measure $m_G$. 
An increasing family of bounded Borel subsets $\cU_r$ ($r\in \RR_+$ or $r\in \NN_+$) 
of $G$ will be called
{\it coarsely admissible} if 
\begin{itemize}
\item For every 
bounded $B\subset G$, there exists $c=c(B)>0$ such that for all $r\ge 1$,
\begin{equation}\label{eq:B_max}
B\cdot \cU_r\cdot B \subset \cU_{r+c}.
\end{equation}

\item For every $c>0$, there exists $D=D(c)>0$ such that for all $r \ge 1$,
\begin{equation}\label{eq:B_max2}
m_G(\cU_{r+c})\le D\cdot  m_G(\cU_{r}).
\end{equation}
\end{itemize}
}
\end{definition}
Let us note that for a length function  $\tilde{\cG}$ on $G$ subject to the standing assumptions stated above, the balls $\tilde{B}_r$ constitute a family of sets satisfying (\ref{eq:B_max}). Indeed, for any given bounded $B$, we have $B\subset \tilde{B}_{s} $ for some $s \ge 1$, and the equation 
$$\tilde{\cG}(b_1gb_2)\le \tilde{\cG}(b_1)+\tilde{\cG}(g)+ \tilde{\cG}(b_2)$$ 
implies when $b_1,b_2$ range over  $\tilde{B}_{s}$ that 
$$ B\tilde{B}_{r}B\subset \tilde{B}_{s} \tilde{B}_{r}\tilde{B}_{s}\subset \tilde{B}_{r+2s}\,.$$

Furthermore, if $\tilde{B}_{r}$ has almost exact polynomial-exponent volume growth, then (\ref{eq:B_max2}) is clearly satisfied as well. 

Note that if $\tilde{B}_{r}$ has almost exact polynomial-exponential volume growth, then the family of spherical shells $\tilde{S}\tilde{S}_r$, for any fixed width $L\ge L_0$ (for some fixed positive $L_0$) has the following property. 
Since $\tilde{B}_{r}$ satisfies (\ref{eq:B_max}) then (since $e\in  \tilde{B}_{s} $) : 
$$\tilde{B}_{r-2s}\subset \tilde{B}_{s} \tilde{B}_{r-2s}\tilde{B}_{s} \subset  \tilde{B}_{r} \subset \tilde{B}_{s} \tilde{B}_{r}\tilde{B}_{s}\subset \tilde{B}_{r+2s}\,,$$
and so 
\begin{equation}\label{coarse-shells}\tilde{S}\tilde{S}_r \subset  \tilde{B}_{s}\cdot (\tilde{S}\tilde{S}_r )\cdot \tilde{B}_{s}= \tilde{B}_{s}  \left( \tilde{B}_{r+L} \setminus  \tilde{B}_{r} \right) \tilde{B}_{s} \subset  \left(\tilde{B}_{r+L+2s}\setminus  \tilde{B}_{r-2s}\right)\,.
\end{equation} 
Finally, we will assume that $\tilde{S}\tilde{S}_r$ satisfies the rapid decay property for convolution norms, namely for $r > 1$
\begin{equation}\label{RRD-G}  
\norm{\lambda_G\left(\chi_{\tilde{S}\tilde{S}_r}\right)}_{L^2(G)\to L^2(G)} \le C r^{\bf b} (m_G(\tilde{S}\tilde{S}_r))^{1/2}\,.
\end{equation} 

We now turn to show that almost exact solution to the lattice point counting problem in a coarsely admissible family of sets whose associated averages have rapid decay is sufficient  the deduce the comparison arguments that we will use. The following result generalizes the discussion of \cite[\S 3.3, Cor. 3.6]{N98}.  

\begin{thm}\label{RRD-lattice-sbgps}
Let $\tilde{\cG}$ be a length function on $G$ satisfying our standing assumptions. Assume that $\tilde{\cG}$ defines balls  $\tilde{B}_r$ with almost exact polynomial-exponential volume growth, and that the lattice $\Gamma$ in $G$ has an almost exact solution to the lattice point counting problem in $\tilde{B}_r$.  Assume also the property of rapid decay with polynomial parameter ${\bf b}$ for the spherical shells $\tilde{S}\tilde{S}_r$ on $G$ of width at least a fixed $ \tilde{L}\ge \tilde{L}_0> 0$, w.r.t. to the length function $\tilde{\cG}$.  Then $\Gamma$ satisfies the property of rapid decay  of spherical shells $SS_r$ of width at least a fixed $L \ge L_0> 0$, with polynomial parameter $\bf b$. 
\end{thm}
\begin{proof}
 Let $F$ and $\psi$ be continuous compactly supported function on $G$. Consider the $\Gamma$-periodization $\phi$ of $\psi$, given by 
$\phi(x\Gamma)=\sum_{\gamma\in \Gamma} \psi(x\gamma)$, which satisfies $\phi\in L^2(G/\Gamma)$. Then the following operator norm inequality holds:
\begin{equation}\label{induction} 
\norm{\lambda_\Gamma\left((\psi^\ast\ast F\ast \psi)|_\Gamma\right)}_{\ell^2(\Gamma)\to \ell^2(\Gamma)}\le 
\norm{\phi}^2_{L^2(G/\Gamma)} \norm{\lambda_G(F)}_{L^2(G)\to L^2(G)}
\end{equation} 
 This inequality is a consequence of the inequality stated in \cite[\S 7]{BS93} (see the discussion following Lemma 7.3)  and described there as the dual version of  \cite[Prop. 1.2]{CS96}, see also the discussion in \cite[\S 3.3, Cor. 3.6]{N98}.

Fix $r_1 > 1$,  so that $m_G(\tilde{B}_{r_1})>0$ and let $u_1=\chi_{\tilde{B}_{r_1}}$ 
and $\psi=u_1\ast u_1$. Then $\psi$ is continuous function with $\psi(x) \ge C^{-1}> 0$ on $\tilde{B}_{r_0}$ for some $C^{-1}> 0$ and some $0 <  r_0 \le r_1$ such that $m_G(\tilde{B}_{r_0}) > 0$. Therefore $u_0=\chi_{\tilde{B}_{r_0}}\le C\psi$ on $G$, and for any non-negative bounded compactly supported function $F$ on $G$,  the  inequality $u_0^\ast \ast F\ast u_0\le C^2\psi^\ast\ast F\ast \psi $ holds at every point in $G$. 
As a result, restricting to $\Gamma$ we have 
\begin{equation}\label{conv}
\left(u_0^\ast\ast F\ast u_0 \right)|_\Gamma \le C^2\left(\psi^\ast\ast F\ast \psi\right)|_\Gamma\,.
\end{equation}

Note that in general, if $0\le h_1\le  h_2$ for (say) finitely supported non-negative functions on $\Gamma$, then $\norm{h_1\ast}_{\ell^2\to \ell^2} \le  \norm{h_2\ast}_{\ell^2\to \ell^2}$. Therefore we conclude that 
\begin{equation}\label{norm-comp}
\norm{\lambda_\Gamma\left(\left(u_0^\ast\ast F\ast u_0 \right)|_\Gamma\right)}_{\ell^2\to \ell^2}  \le  C^2\norm{\lambda_\Gamma\left(\left(\psi^\ast\ast F\ast \psi\right)|_\Gamma}\right)_{\ell^2\to \ell^2}\,.
\end{equation} 
By (\ref{induction}), the r.h.s of (\ref{norm-comp}) is bounded by $C_1 \norm{\lambda_G(F)}_{L^2(G)\to L^2(G)}$.  

We now consider two families of sets, as follows.  
The first is the family of spherical shells in $\Gamma$ of width $L$ consisting of lattice points satisfying 
$$SS_r=\set{\gamma \in \Gamma\,;\, r\le \cG(\gamma)< r+L}. $$
The second is the family of spherical shells of width $L+4r_0$ in $G$ given by (for $r > 3r_0$)  
$$\tilde{S}\tilde{S}^{\prime}_r=\set{g\in G\,;\, r-2r_0\le \cG(g)< r+L+2r_0}.$$

We proceed to bound the operator defined on $\ell^2(\Gamma)$ by the discrete sum over $SS^\prime_{r}=\Gamma\cap \tilde{S}\tilde{S}^{\prime}_r$, namely the convolution operator defined by the characteristic function $h_{1,r}$ of the set $SS^\prime_{r}$. For a bounded subset $U\subset G$ of positive measure 

$$u_0^\ast\ast \chi_U\ast u_0(g)=
\int_{x\in G}\int_{y\in G} 
\chi_{\tilde{B}_{r_0}}(x^{-1})\chi_U(xgy^{-1})\chi_{\tilde{B}_{r_0}}(y)dm_G(x)dm_G(y)
$$
Now let $\tilde{U}= \tilde{S}\tilde{S}^\prime_r$, and recall that $\tilde{B}_r$ is a coarsely admissible family. Using (\ref{coarse-shells}) with $s=r_0$ 
we conclude that  
$$\tilde{B}_{r_0} \cdot\tilde{S}\tilde{S}_{r} \cdot \tilde{B}_{r_0} \subset 
 \left(\tilde{B}_{r+L+2r_0}\setminus  \tilde{B}_{r-2r_0}\right)= \tilde{S}\tilde{S}^\prime_{r}\,.$$
It follows that for $\gamma \in SS_{r}$ and $x,y\in \tilde{B}_{r_0}$, we have 
$\chi_{\tilde{S}\tilde{S}^\prime_{r}}(x\gamma y^{-1})=1$, and hence for $\gamma \in \tilde{S}\tilde{S}_{r}$, we have $u_0^\ast\ast \chi_{\tilde{U}}\ast u_0(\gamma)\ge m_G(\tilde{B}_{r_0})^2 > 0$. 
We conclude that 
\begin{equation}\label{restr-estim} 
\chi_{SS_{r}}=\chi_{\tilde{S}\tilde{S}_{r}}|_\Gamma\le C_2\left(u_0^\ast \ast 
\chi_{\tilde{S}\tilde{S}^\prime_r}\ast u_0\right)|_\Gamma\,.
\end{equation}

It follows from (\ref{conv}) that  
$$h_{1,r}=\chi_{\tilde{S}\tilde{S}_{r}}|_\Gamma\le C_3(\psi\ast \chi_{\tilde{S}\tilde{S}^\prime_r}\ast \psi)|_\Gamma=h_{2,r}$$
 and hence that the operator norm of the convolution operator $h_{1,r}\ast$ in $\ell^2(\Gamma)$ is dominated by the operator norm of the convolution operator $h_{2,r}\ast$. 

The operators $(\psi\ast\chi_{\tilde{S}\tilde{S}^\prime_r}\ast \psi)|_\Gamma$ have operator  norms in the regular representation of $\Gamma$ which satisfy the estimate (\ref{induction}). By assumption, the operator norm of 
$\lambda_G(\chi_{\tilde{S}\tilde{S}^\prime_r})$ on $L^2(G)$ appearing in the r.h.s. of (\ref{induction}) satisfy the property of rapid decay in $L^2(G)$ stated in (\ref{RRD-G}). Normalizing the operators we conclude that for $r \ge 1$
$$\norm{\lambda_G\left(\frac{\chi_{\tilde{S}\tilde{S}^\prime_r}}{m_G(\chi_{\tilde{S}\tilde{S}^\prime_r})}\right)}_{L^2\to L^2}\le C_4 r^{\bf b } m_G(\tilde{S}\tilde{S}^\prime_r)^{-1/2}\le C_5 r^{{\bf b}-\frac{d}{2}}q^{-r/2}\,.
$$ 
We now normalize the convolution operators $h_{1,r}$ supported on $SS_{r}$ and use our assumptions of almost exact polynomial-exponential growth of spherical shells on $G$ and hence also on $\Gamma$ with the same parameters $(d,q)$. This uses the almost exact solution of the lattice point counting problem as noted in Lemma \ref{restricted-gauge-volume}. Finally we conclude that for $r$ at least a fixed positive constant  
\begin{equation}\label{RRD-Gamma}
\norm{\lambda_G\left(\frac{\chi_{SS_r}}{\abs{SS_r}}\right)}_{\ell^2(\Gamma)\to \ell^2(\Gamma)}\le C_6 r^{{\bf b}-\frac{d}{2}}q^{-r/2} \le C_7 r^{\bf b } \abs{ SS_r}^{-1/2}\,.
\end{equation}
This concludes the proof of Theorem \ref{RRD-lattice-sbgps}. 
\end{proof} 
\subsection{Additional remarks and discussion of radial rapid decay} 

Let us note the following relevant prior results and references.

\begin{enumerate}
\item M. Perrone considered the balls $\tilde{B}_n$ defined by the (integer-valued) word metric associated with certain bounded symmetric generating set in a connected semisimple Lie group $G$ with finite center. For an irreducible uniform lattice, rapid decay for the averages on  $\Gamma \cap \tilde{B}_r$ is established in 
\cite[Cor. 1.4]{Pe09}.  Rapid decay for Riemannian spherical shells $\tilde{S}\tilde{S}_{n\epsilon} $ intersected with the lattice is established  in \cite[Cor. 3.8]{Pe09}. 
 
\item A. Valette  considered a group acting simply transitively on the vertices of an $\tilde{A}_n$-building. Rapid decay for the sphere (and general radial) averages, with respect to the 
combinatorial distance function on the $1$-skeleton of the building is established in \cite[Thm. 1]{Va97}. Rapid decay for general uniform lattices in simple algebraic groups over a non-Archimedean local field, w.r.t. the combinatorial distance on the Bruhat-Tits building is stated in \cite[Thm. 1.5]{Pe09}. 
\item  Boyer, Lobos and Pittet \cite{BLP22} have recently established a general dynamical criterion to establish radial rapid decay, using the action of a discrete group on its boundary. It is shown that with respect to a suitable length function on irreducible lattices of semisimple Lie groups, radial rapid decay holds. The argument uses the lattice action on the Poisson boundary $G/P$, estimates of the Harish Chandra $\Xi_G$-function, and operator norm bounds in the associated unitary representation on the boundary, which was established earlier in \cite{BLP17}.  The method is also applied to a certain non-uniform lattice $\Lambda$ in the automorphism group of a products of a regular tree with itself, in order to show that radial rapid decay for a length function quasi-isometric to a word metric on $\Lambda$ does hold, but that $\Lambda$ does not have the rapid decay property, because it contains a solvable group with exponential growth (see \cite{Ch17}).   

\item A wealth of additional examples which satisfy rapid decay of spherical shells but not rapid decay is provided by Theorem \ref{RRD-lattice-sbgps} together with Theorem \ref{lattice-ss-gps} below. To demonstrate this fact, consider the (non-uniform) lattices $SL_n(\ZZ)$, $n \ge 3$, and the length function arising from restriction of the Riemannian distance on the symmetric space, discussed below. Here Theorem \ref{RRD-lattice-sbgps} asserts that Riemannian spherical shells restricted  to $SL_n(\ZZ)$ satisfy rapid decay, and the same holds for the balls. But these lattices do not satisfy property RD, as they contain a solvable group with exponential growth.  Finally, the Riemannian distance restricted to $SL_n(\ZZ)$ is in fact quasi-isometric to a word length on the lattice, as shown by \cite{LMR00}. 
\item A discussion of some consequences of radial rapid decay, particularly for word-hyperbolic groups, in the context of positive unitary representation with a spectral gap appears in Boucher \cite[\S 5]{B24}. 
\end{enumerate}

\subsection{Weak-type maximal inequality for lattice subgroups}
We now turn to the proof of the following.

\begin{thm}\label{lattice-ss-gps} 
Let $G$ be a connected semisimple Lie group with finite center, $\Gamma$ any lattice subgroup. Fix an origin $o$ in the symmetric $G/K$ 
and define the length function 
$\cG(\gamma)=D(\gamma o,o)$, where $D$ is the $G$-invariant Riemannian metric arising from the Killing form. Then the Hardy-Littlewood maximal function on $\Gamma$ associated with this length function satisfies a maximal inequality in $\cL(\log \cL)^{\bf c}(\Gamma)$. Here ${\bf c}$ is an explicit function of the root system of 
$G$ only, given explicitly below. 
\end{thm} 
\begin{proof}
We verify that the assumptions of Theorem \ref{RRD-lattice-sbgps} are satisfied, and then the desired conclusion follows from Theorem \ref{main}. 

I) The Riemannian balls in symmetric space $G/K$ defined by the Killing form have almost exact polynomial-exponential growth, and therefore the same holds for the balls $\tilde{B}_r\subset G$ defined by the admissible length function $\tilde{\cG}(g)=D(go,o)$, whose restriction to $\Gamma$ is $\cG$. In fact, the volume growth is {\it exact}, and  the asymptotic formula for the volume of Riemannian balls in symmetric space  is classical. We refer to \cite[\S 5]{GO07} for a succinct account, and to 
\cite[Thm. 6.2]{Kn95} for the even stronger fact that the Riemannian area of the Riemannian spheres also has {\it exact} polynomial-exponential growth.

II) The lattice point counting problem in Riemannian balls in symmetric space has in fact an {\it effective} solution, as can be easily be deduced from the method developed in \cite{DRS93} (which was applied there to the case of certain linear norms). For a different and more general proof, see \cite{gn_counting}. We need here only the fact that the solution is almost exact, which of course follows from the fact that it is exact established e.g. in \cite[Thm. 1.2]{DRS93}, or  \cite[Thm. 1.4]{EM93}. Indeed the symmetric space is affine symmetric w.r.t. the compact subgroup $K$, and the Riemannian balls are well-rounded in the sense of \cite[Prop. 1.3]{EM93}.  

III) Since statement  I) implies the almost exact polynomial-exponential growth of the balls $\tilde{B}_r$ in $G$, and statement II) implies the almost exact 
solution to the lattice point counting in them, Lemma \ref{restricted-gauge-volume}
implies that the spherical shells $SS_r$ in $\Gamma$ have almost exact polynomial-exponential growth with the same parameters (for any given $L > L_0$ defined in the Lemma). 

IV) Since the length function $\tilde{\cG}$ defining the balls $\tilde{B}_r$, namely $\tilde{\cG}(g)=\tilde{D}(g,e)$  is given by an invariant metric, we have  $\tilde{\cG}(b_1 gb_2)\le \tilde{\cG}(b_1)+\tilde{\cG}(g)+\tilde{\cG}(b_1)$. Since in addition the volume growth of $\tilde{B}_r$ is almost exactly polynomial-exponential, it follows that $\tilde{B}_r$ are coarsely admissible, see Def. (\ref{coarse-def}). 

V) To conclude the verification that the assumptions of Theorem \ref{RRD-lattice-sbgps} are satisfied, we must show that spherical shells $\tilde{S}\tilde{S}_r$ in $G$
(of any width $L\ge  L_1 > 0$) satisfy rapid decay of convolution norms. 

To begin with, we note that in the present case the norm of a bi-$K$-invariant average viewed as a convolution operator on $L^2(G)$ is simply given by the integral of the Harish Chandra $\Xi$-function of $G$ w.r.t. this average, namely for any bounded $K$-bi-invariant set $\tilde{U}$ of positive measure 
\begin{equation}\label{spectral norm}
\norm{\lambda_G(1_{\tilde{U}})}_{L^2(G)\to L^2(G)}=\int_{\tilde{U}}\Xi_G(g)dm_G(g)
\le m_G(\tilde{U})^{1/2} \left(\int_{\tilde{U}}\Xi_G(g)^2dm_G\right)^{1/2}
\end{equation}  
For a succinct account of this classical fact see \cite[Prop. 4.3, Thms. 4.3, 4.5]{CPSC}. 
 
 We will now use the notation of \cite[\S 2.2]{GV88} for the closed positive Weyl chamber  $\overline{\mathfrak{a}+}$, its exponential $\exp\overline{\mathfrak{a}+}=\overline{A_+}$, and the Cartan polar coordinates decomposition $G=K\overline{A_+}K$, namley $g=k_1 \cdot \exp(H(g))\cdot k_2$, the radial part being denoted $g\mapsto H(g)$. The expression of Haar measure in polar coordinates is given by 
  (see \cite[Prop. 2.4.6]{GV88})  
 \begin{equation}\label{int-formula} \int_G f(g)dm_G(g)=\int_K \int_{\overline{\mathfrak{a}+}}\int_K f(k_1\cdot \exp H \cdot  k_2)J(H)
dk_1 dH dk_2\end{equation} 
 with the density given by $J(H)=\prod_{\alpha\in \Delta^+}\left(\sinh\alpha(H)\right)^{m_\alpha}$, 
 $\Delta^+$ the set of positive roots of $\mathfrak{a}$ in $\mathfrak{g}$, and $m_\alpha$ the multiplicity of the root $\alpha$.

 Recall the well-known Harish Chandra estimates of the Harish Chandra $\Xi$-function (see e.g. \cite[Thms. 4.6.4, 4.6.5]{GV88}, and more succinctly their refinement by Anker \cite{An87}), given by, for all $H\in \overline{\mathfrak{a}+}$ :
 \begin{equation}\label{HC-estimate}  e^{-\rho(H)}\le \Xi_G(\exp H) \le C e^{-\rho(H)}\left(1+\norm{H}\right)^{\abs{\Delta^+_0}} \end{equation} 
where $\rho$ is half the sum of the positive roots (with multiplicities), $\Delta^+_0$ the set of short positive roots, and $\norm{H}$ is the Cartan-Killing norm on $\mathfrak{a}$. 

Finally recall also that the Riemannian distance $D$ on the symmetric space $G/K$ satisfies $D(go, o)=D(k_1 \cdot \exp(H(g))\cdot k_2 \cdot o, o)=\norm{H}$ where $o$ is the coset $K\in G/K$. The Riemannian ball of radius $r$ in $G/K$ lifted to $G$ is therefore given by $\tilde{B}_r=\set{g\in G\,;\, \norm{H(g)}\le r}$.  
 
 Since $\Xi_G$ is a $K$-bi-invariant function, it follows from the integration formula (\ref{int-formula}) and the estimate of the Harish Chandra function (\ref{HC-estimate}))
 that for $r \ge 1$ 
 \begin{equation}\label{L2+epsilon} 
 \int_{\tilde{B}_r} \Xi_G^2(g) dm_G \le C r^{\dim \mathfrak{a}+2\abs{\Delta^+_0}}
 \end{equation}
 Indeed, 
since $J(\exp(H)\le \exp(2\rho(H)$ for $H\in \overline{\mathfrak{a}^+}$ (see e.g. the proof of \cite[Prop. 4.6.12]{GV88})
 $$ \int_{\tilde{B}_r} \Xi_G^2 (g) dm_G=
 \int_{H\in\overline{\mathfrak{a}+}\,;\, \norm{H}\le r} \Xi_G^2 (\exp(H)) J(H)dH$$   
 $$\le  \int_{H\in\overline{\mathfrak{a}+}\,;\, \norm{H}\le r}\left(1+\norm{H}\right)^{2\abs{\Delta^+_0}}dH\le C_1 \int_0^r t^{\dim \mathfrak{a}-1} (1+t)^{2\abs{\Delta^+_0}}dt\le C r^{\dim \mathfrak{a}+2\abs{\Delta^+_0}}$$

 It follows using (\ref{spectral norm}) that the ball averages $\tilde{B}_r$, $r \ge 1$ satisfy the norm estimate 
\begin{equation}\label{Xi-function} 
\norm{\lambda_G(\tilde{\beta}_r)}_{L^2\to L^2}\le m_G(\tilde{B})^{-1/2} \left(\int_{\tilde{B}_r}\Xi_G(g)^2dm_G\right)^{1/2}\le C_G r^{{\bf b}} m_G(\tilde{B}_r)^{-1/2}\,,
\end{equation} 
where ${\bf b}={\bf b}_G=\frac12 \dim \mathfrak{a}+\abs{\Delta^+_0}$. 

VI) Finally, the fact that the convolution averages on $\tilde{S}\tilde{S}_r$ satisfy the rapid decay estimate with the same polynomial bound easily follows from (\ref{spectral norm}) with $\tilde{U}=\tilde{S}\tilde{S}_r$. Indeed, we can then extend the integral on the r.h.s of (\ref{spectral norm}) from $\tilde{U}$ to $\tilde{B}_{r+L}$, and the normalization by the measure of $\tilde{S}\tilde{S}_r$ will give the desired estimate, using the fact that the volume growth of $\tilde{B}_r$ is almost exact. 

This concludes the proof of Theorem \ref{lattice-ss-gps}. 
\end{proof}
\begin{remark}
The conclusion of Theorem  \ref{lattice-ss-gps} holds far beyond 
the cases stated in it. A weak-type maximal inequality in $\cL(\log \cL)^{\bf c}(\Gamma)$ holds for lattices $\Gamma$ in $S$-algebraic groups $G$ over locally compact second countable totally disconnected fields, with respect to suitable distance functions. It also holds for natural families of radial averages on the associated Bruhat-Tits buildings. This can be established without difficulty using the arguments developed in the proof of Theorem  \ref{lattice-ss-gps} for specific examples such as radial averages on $SL_n(\QQ_p)$, its lattice subgroups, and the associated building.  The algebraic and analytic background and details required to establish such a result in full generality require quite elaborate expositions and will have to be taken up elsewhere. 
\end{remark}

\section{The spherical coarse-median inequality}
We now formulate the following inequality which will allow us to prove a maximal inequality for the Hardy-Littlewood operator
using it. 

\subsection{Spherical coarse median inequality}

Consider the following definition for a general countable group $\Gamma$, 
and a left-invariant integer-valued metric $d_\cG(x,y)=\cG(x^{-1}y)$ on $\Gamma$, with associated spheres $S_r=\set{w\in \Gamma \,;\, \cG(w)=r}$ (namely spherical shells with width $L=1$). 

\begin{definition}[Spherical coarse median inequality]\label{median-structure} 
$\Gamma$ satisfies the {\em spherical coarse median inequality of rank $(d_2+1)$} if the following holds for some 
constant $C>0$. Given any $r\in \NN_0$, $0\le m \le r$, and $i, j\in \NN_0$ satisfying $i=j+r-2m$, any two non-empty sets $E_j\subset S_j$, $F_i\subset S_i$ satisfy 
\begin{equation}\label{median-estimate} 
\abs{\set{(x,y)\in E_j\times F_i\,;\, d_\cG(x,y)=r}}\le  C_0 r^{d_2}\min\set{\abs{S_{r-m}}\abs{E_j}, \abs{S_m}\abs{F_i}}
\end{equation}
\end{definition} 

Note that (\ref{median-estimate}) generalizes the case where $d_2=0$ which is utilized in the proof of \cite[Lemma 5.1]{NT09}. Some comments explaining this terminology  appear at the end of \S 8.

\subsection{Spherical coarse median inequality and rapid-decay of spherical correlations  }

In the present section we will consider word metrics (and more generally, integer-valued length functions) and prove the following result.

\begin{thm}[Main theorem II : Spherical coarse median inequality]\label{thm-HL-d>0}
Let $(\Gamma,\cG)$ satisfy the spherical coarse median property with parameter $d_2$  and have  almost exact polynomial-exponential  growth of spheres with parameters $d$ and $q$. Then the Hardy-Littlewood operator satisfies the weak-type  $\cL\left(\log \cL\right)^{\bf c}(\Gamma)$ maximal inequality, with ${\bf c}=2d_2+d$.
In particular, when $d_2=0=d$ the Hardy-Littlewood operator satisfies the weak-type $(1,1)$-maximal inequality in $\ell^1(\Gamma)$. 
\end{thm} 

The proof of Theorem \ref{thm-HL-d>0} is based on the following proposition, which  constitutes a direct generalization of the novel radial geometric argument developed in \cite{NT09} for the case of rooted or regular tree.  

 \begin{prop}\label{convolution-ineq} 
 If $(\Gamma,\cG)$ satisfies the spherical coarse median inequality with parameter $d_2$, and has  almost exact polynomial-exponential growth with parameters $(q,d)$, then $\Gamma$ satisfies the property of rapid-decay of spherical correlations (\ref{weak-RD}), with polynomial parameter  ${\bf b}=d_2+\frac12 d$. 
\end{prop}
\begin{proof}
Given two non-empty finite sets $E$ and $F$, we set $E_j=E\cap S_j$ and $F_j=F\cap S_j$. Now write 
$$
\abs{\set{u\in E, v\in F\,;\, d_\cG(u,v)=r}}=
$$
$$
\sum_{m=0}^r\sum_{i,j\in \NN_0\,;\, i=j+r-2m}
\abs{\set{(x,y)\in E_j\times F_i\,;\, d_\cG (x,y)=r}}.
 $$

Using the assumption that $(\Gamma, \cG)$ satisfies (\ref{median-estimate}) 
it suffices to prove 
$$
\sum_{m=0}^r\sum_{i,j\in \NN_0\,;\, i=j+r-2m}
r^{d_2}\min\set{\abs{S_{r-m}}\abs{E_j}, \abs{S_m}\abs{F_i}}$$
$$
 \le C r^{d_2+\frac12 d} \sqrt{\abs{E}}\sqrt{\abs{F}}\sqrt{\abs{S_r}}
$$

Set $u_j=\frac{\abs{E_j}}{\abs{S_j}}$, $v_j=\frac{\abs{F_j}}{\abs{S_j}}$, for $j \ge 0$ and $u_j=v_j=0$ for $j < 0$.

Then 
$$\sum_{j\in \NN_0} \abs{S_j}u_j=\abs{E} \text{  and  } \sum_{j\in \NN_0} \abs{S_j}v_j=\abs{F},$$
and 
$$\sum_{m=0}^r\sum_{i,j\in \NN_0\,;\, i=j+r-2m}
r^{d_2}\min\set{\abs{S_{r-m}}\abs{E_j}, \abs{S_m}\abs{F_i}}\le 
$$
$$\le C_1\sum_{m=0}^r \sum_{i,j\in \NN_0\,;\, i=j+r-2m}
r^{d_2}\min \set{(r-m)^{d} q^{r-m}\abs{E_j}, m^{d} q^m\abs{F_i}}\le 
$$
$$\le C_2\sum_{m=0}^r\sum_{i,j\in \NN_0\,;\, i=j+r-2m}
r^{d_2}\min \set{(r-m)^{d} q^{r-m}j^{d} q^j u_j, m^{d} q^m i^{d} q^i v_i}
$$
Noting that $m+i=r-m+j$, we can factor out $q^{m+i}$, and using 
$\frac12 (r+i+j)=m+i$, the last expression equals 
\begin{equation}\label{intermediate} C_2 r^{d_2}q^{r/2}\sum_{m=0}^r\sum_{i,j\in \NN_0\,;\, i=j+r-2m}q^{(i+j)/2}
\min \set{\left((r-m)j\right)^{d}  u_j, \left(mi\right)^{d} v_i}.
\end{equation} 

We now claim that  the double sum appearing in the last expression satisfies the following estimate : 

$$\sum_{m=0}^r\sum_{i,j\in \NN_0\,;\, i=j+r-2m}q^{(i+j)/2}
\min \set{\left((r-m)j\right)^{d}  u_j, \left(mi\right)^{d} v_i}\le 
$$
$$\le 2 r^{d}
\left(\sum_{j\in \NN_0} j^{d} q^j u_j \right)^{1/2}\left(\sum_{i\in \NN_0}
 i^{d} q^i v_i\right)^{1/2}\le 2 C_3 r^{d}\sqrt{\abs{E}}\sqrt{\abs{F}}
$$

To prove the claim, note that for a fixed $r$ and fixed $m$ with $0\le m \le r$, the foregoing inner sum ranges over
$i,j\in \NN_0$ with $i=j+(r-2m)$, which is a line in $\NN_0\times \NN_0$ parallel to 
$i=j$. As $m$ varies from $0$ to $r$, these lines are disjoint. Therefore 
since $0\le m\le r$  
$$\sum_{m=0}^r\sum_{i,j\in \NN_0\,;\, i=j+r-2m}q^{(i+j)/2}
\min \set{\left((r-m)j\right)^{d}  u_j, \left(mi \right)^{d} v_i}\le 
$$

$$\le r^{d}\sum_{m=0}^r\sum_{i,j\in \NN_0\,;\, i=j+r-2m}q^{(i+j)/2}
\min \set{j^{d}  u_j, i^{d} v_i}\le 
$$
$$
\le r^{d} \sum_{i,j \in \NN_0} q^{(i+j)/2}
\min \set{j^{d}  u_j, i^{d} v_i}.  
$$
Since the minimum chooses either $i$ or $j$, for any $\alpha \in \RR$  the last sum is bounded by 
$$\le  r^{d}\left( \sum_{i,j \in \NN_0\,,\, i < j+\alpha} q^{(i+j)/2} j^{d}  u_j+ \sum_{i,j \in \NN_0\,,\, i \ge j+\alpha}  q^{(i+j)/2}  i^{d} v_i \right),  
$$
and using the two conditions $i < j+\alpha$ and $i \ge j+\alpha$ we obtain the bound 
$$\le  r^{d}\left(\sum_{j\in \NN_0} q^{j+\frac{\alpha}{2}}  j^{d} u_j +
 \sum_{i\in \NN_0} q^{i-\frac{\alpha}{2}}  i^{d} v_i\right).
$$
By our assumption on sphere growth, 
$k^{d}q^{k}\frac{\abs{E_k}}{\abs{S_k}}\le C_4 \abs{E_k}$ and therefore the latter expression satisfies  
$$ =r^{d}\left(q^{\frac{\alpha}{2}} \sum_{j\in \NN_0} q^{j}  j^{d} u_j +
 q^{-\frac{\alpha}{2}}\sum_{i\in \NN_0} q^{i}  i^{d} v_i\right)\le  
C_5 r^{d}\left( q^{\frac{\alpha}{2}} \abs{E}+q^{-\frac{\alpha}{2}}\abs{F}\right)
$$
The minimum of the expression in brackets is obtained by differentiation w.r.t. $\alpha$. The minimum is obtained when $q^\alpha=\frac{\abs{F}}{\abs{E}} $, and its value 
is then $2\sqrt{\abs{E}}\sqrt{\abs{F}}$. 

Finally, we combine (\ref{intermediate}) with the estimate just established. Since $r^{d/2}q^{r/2}\le \sqrt{C_6}\sqrt{\abs{S_r}}$ we conclude that the bound established is 
$$C_7 r^{d_2} r^{d}q^{r/2} \sqrt{\abs{E}}\sqrt{\abs{F}}\le C_8 r^{d_2+\frac12 d} \sqrt{\abs{E}}\sqrt{\abs{F}}\sqrt{\abs{S_r}}. 
$$
This completes the proof of  Proposition  \ref{convolution-ineq}. 
\end{proof}

{\it Proof of Theorem \ref{thm-HL-d>0}}. Under the assumptions of Theorem \ref{thm-HL-d>0}, Proposition \ref{convolution-ineq} implies that  the sphere averages satisfy rapid decay of spherical correlations with polynomial parameter ${\bf b}=d_2+\frac12 d $.  By Theorem \ref{main}, since we are assuming that the growth of the spheres is almost exact polynomial-exponential it follows that the Hardy-Littlewood operator satisfies the maximal inequality stated.  
\qed

\section{Hyperbolic groups}

We apply Theorem \ref{thm-HL-d>0} to hyperbolic groups.
We 
first prove the coarse median inequality of rank $1$ (namely $d_2=0$) for hyperbolic groups. 
This result generalizes the case of free groups with a free generating sets which was  established by \cite[inequality (76)]{NT09}, and the proof follows the same idea. 

\begin{prop}\label{Coarse Med Hyp}  Non-elementary hyperbolic groups satisfy the spherical coarse median inequality of rank $1$ (namely with parameter $d_2=0$), for every word metric associated with a  
finite symmetric generating set. 
\end{prop}
\begin{proof}
Let $\Gamma$ be a non-elementary hyperbolic group with a finite symmetric generating set $S$, and denote by the Cayley graph of $\Gamma$ with respect 
to the associated length function $\cG=\cG_S=\abs{\cdot}_S$ by $(\Gamma, \cG)$. This is $\delta$-hyperbolic graph
for a constant $\delta\ge0$, which means that every geodesic
triangle is $\delta$-thin in the sense that each side is contained
in the $\delta$-neighbourhood of the union of the other two sides. 
For $x,y \in \Gamma$, we denote a choice of geodesic  between 
them by $[x,y]$. As usual, denote the distance between $x,y \in \Gamma$ by $d_{\cG}(x,y)=\cG(x^{-1}y)$.

Let $e$ denote the identity element of $G$, and represent  $\Gamma=\cup_{r\ge 0} S_r$, the disjoint union of the spheres
$S_r=\{g \in \Gamma \,;\, \cG(g)=r \}$.
A growth estimate for the spheres was established  in \cite{Co93}, as follows.
There exist $C\ge 1$ and $q>1$ such that 
for any $r$ we have 
\begin{equation}\label{ineq.Co}
\frac{1}{C}q^r \le\abs{S_r} \le C q^r.
\end{equation}
Therefore $G$ has almost exact  exponential growth
with parameters $q$ and $d=0$.

\begin{figure}[htbp]
\hspace*{-3.3cm}     
\begin{center}
                                                      
   \includegraphics[scale=0.5]{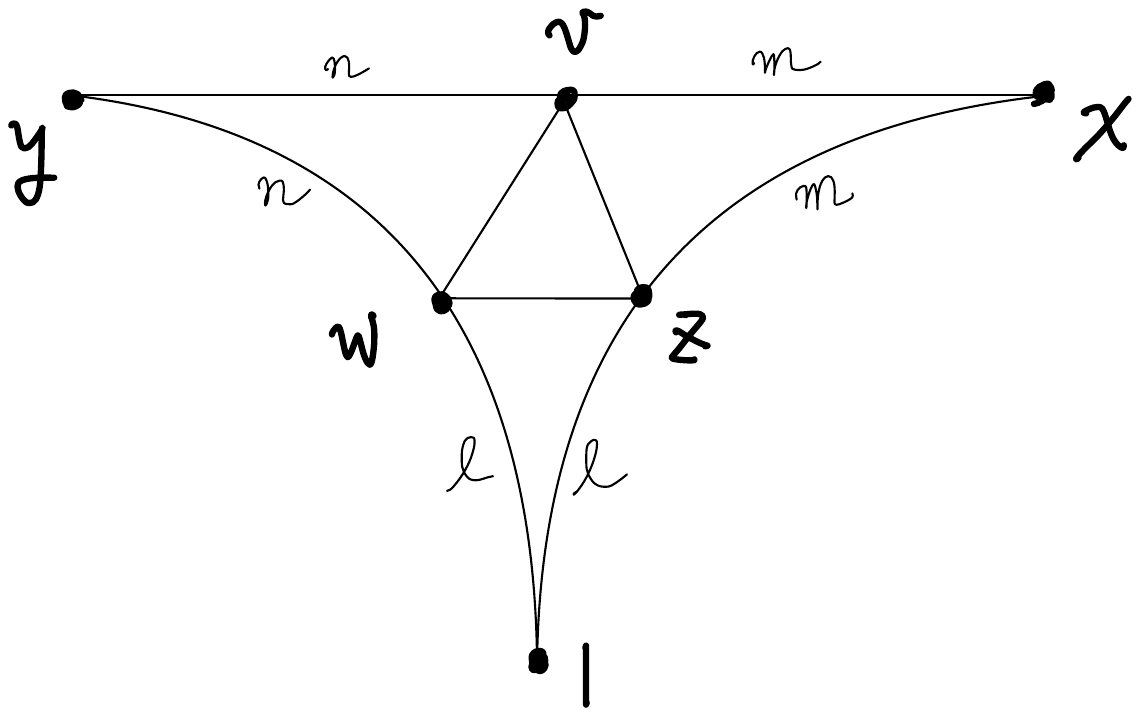}
   \end{center}
\caption{ $\delta$-thin triangle. $1$ is the identity element $e$.}
  \label{triangle}
\end{figure}

{\bf Spherical Coarse Median Inequality} : Let  $E_j \subset  S_j, F_j \subset S_j$.

Given $i,j \ge 0$ and $r$, define 
$m$ such that $i=j+r-2m$, namely,
$$m=\frac{j+r-i}{2}.$$
By the triangle inequality, we have $0 \le m \le r$.
Here, $m$ could be a half integer, but for simplicity, we pretend
that $m$ is an integer in the following argument. We need a minor
modification in the case that $m$ is an half integer, which we omit.

Set $B=q^{10\delta}$. 
Given $x \in E_j$
there are at most $C B q^{r-m}$ possibilities
for $y \in F_i$ with $r=d(x,y)$.
To see it, set $n=r-m, \ell=i-n=j-m$.
Choose $z\in [1,x], w\in [y,1], v\in [y,x]$ 
such that 
$$d_\cG(y,v)=d_\cG(y,w)=n, d_\cG(v,x)=d_\cG(z,x)=m, d_\cG(w,1)=d_\cG(z,1)=m.$$
See the Figure \ref{triangle}.
Then by $\delta$-hyperbolicity, we have
$$d_\cG(w,z), d_\cG(w,v), d_\cG(v,z) \le 10 \delta.$$
Since $d_\cG(y,w)=n=r-m$, we have 
$$r-m-10\delta \le d_\cG(y,z) \le r-m+10 \delta$$
So, there are at most $C q^{r-m+10\delta}$ possibilities for $y$
by (\ref{ineq.Co}).
Since $B=q^{10\delta}$, we are done.

A similar argument applies to a point $y \in F_i$, 
therefore we have 
\begin{equation}\label{mth} \left| \{ (x,y) \in E_j \times F_{i}: d_\cG(x,y) = r
\}\right| \leq BC \min\left\{ q^{r-m} |E_j|, q^m |F_{i}|
\right\}.\end{equation} 
By (\ref{ineq.Co}), the spherical coarse median inequality holds with polynomial parameter $d_2=0$. 
\end{proof}

We point out that a minor modification of the argument proves that
non-elementary hyperbolic groups satisfy
the spherical coarse median inequality of rank $1$ for any hyperbolic left-invariant integer-valued metric. However, we do not know if these satisfy almost exact exponential growth, namely satisfy (\ref{ineq.Co}).

\begin{thm}[Hyperbolic groups]\label{thm.hyp}
Let $\Gamma$ be a non-elementary hyperbolic group and $S$ any symmetric 
finite generating set, and let $\cG_S=\abs{\cdot}_S$ denote the length-function associated to the word metric. 

Then, the Hardy-Littlewood operator associated with the balls defined by $\cG_S$ satisfies the weak-type $(1,1)$-maximal inequality
in $\ell^1(\Gamma)$. 
\end{thm}

\proof
The Cayley graph for $\Gamma$ w.r.t. $S$ satisfy the coarse median property with parameter $d_2=0$ by 
Proposition \ref{Coarse Med Hyp}.
Also, as noted in the proof, by \cite{Co93} it has almost exact polynomial-exponential growth with parameter $d=0, q>0$.
Now apply Theorem \ref{thm-HL-d>0}.
\qed

\section{Median spaces and the spherical coarse median inequality}

In Proposition \ref{Coarse Med Hyp} we proved that 
a hyperbolic group satisfies the spherical coarse median inequality 
of rank $1$ for any finite generating set. 
In this section, we will establish this inequality for median spaces. 

\subsection{Median space}
We recall the definition of a median space. Let $(X,D)$ be a metric space.
For $x,y \in X$, we say $z \in X$ is between $x,y$ if 
$D(x,z)+D(z,y)=D(x,y)$. Let $C(x,y)$ denote the set
of points between $x,y$. 
$X$ is called a {\em median} space if there is a unique
$m \in C(x,y) \cap C(x,z) \cap C(y,z)$. The point $M=m(x,y,z)$
is called the {\em median point}.

Examples of median spaces include
trees, the $1$-skeleton of a CAT(0) cube complex, $\R^n$ with the $\ell^1$ metrics. 
The $\ell^1$-product of two median spaces
is median. 

We can now state the following result (compare Proposition \ref{Coarse Med Hyp}).

\begin{lem}[Median space and the spherical coarse median inequality] \label{lem-mi}
Let $(\Gamma, D_S)$ be a Cayley graph of a finitely generated group $\Gamma$ w.r.t. the finite symmetric generating set $S$, and the metric $\abs{\cdot}_S=D_S$.  
Assume that it is a median space, and that there exists a polynomial $P(t)$
of degree $\cD$ such that  for any $x,y \in \Gamma$, and $r>0$, 
the set $C(x,y) \cap S(x,r)$ contains at most
$P(r)$ elements, where $S(x,r)=\{z\in \Gamma|
D_S(x,z)=r\}$.

Then $\Gamma$ satisfies the spherical coarse median inequality of rank $(\cD+1)$ (namely with parameter $d_2=\cD$)  
for the metric $D_S$.

\end{lem}
\proof
We adapt the proof of Proposition \ref{Coarse Med Hyp} to this setting.
We use the same notations, such as $n,m,\ell$. 
Since $\Gamma$ is a median space, $m$ is an integer, and the points $z,w,v$ coincide, in a point which is the median
point for $e,x,y$. We denote this point by $M$. 

Let $K$ be a constant such that we have $P(t) \le Kt^{\cD}$ for all $t \ge 0$.
Fix $x \in E_j$, and we consider all $y\in F_i$ with $D_S(x,y)=r$. 
Since $M \in C(x,e)$ with $D_S(x,M)=m \le r$, there are at most
$K r^{\cD} $ possibilities for $M=M(e,x,y)$. For each such $M$, 
there are at most $|S_{r-m}|$ possibilities for $y$ since $D_S(y,M)=r-m$.
Combining them, we obtain
$$ \left| \{ (x,y) \in E_j \times F_{i}: D_S(x,y) = r
\}\right| \le K r^{\cD}  |S_{r-m}||E_j|$$

Similarly, by switching the roles for $x$ and $y$,
 the right hand side is at most $K r^{\cD} |S_m||F_i|$.
 This gives a desired inequality. 
\qed

\subsection{Right-angled Artin groups}
We prove that a right angled Artin group (=RAAG)
satisfies the spherical coarse median inequality. 

\begin{prop}[RAAG satisfies the spherical coarse median inequality]\label{prop.raag}
Let $\Gamma$ be the RAAG defined by a finite graph $L$, 
with the standard generating set $S$.
Let $J$ be the number of vertices in a clique of the largest cardinality in 
the graph $L$ (also known as the dimension of the RAAG). 

Then the Cayley graph of $(\Gamma,S)$ 
satisfies the coarse-median inequality of rank $(J+1)$ (namely with parameter $d_2=J$).

\end{prop}

The proof is based on the following lemma due to \cite{CR05} (p. 334, proof of Theorem 0.4).

\begin{lem}[\cite{CR05}]\label{lemma.CR}
Let $(X,D)$ be a finite dimensional CAT(0) cube complex
of dimension $J$, with the geodesic metric on the 
set of vertices induced from the $1$-skeleton of $X$. 
Then there is a polynomial $P(r)$ of degree $J$
such that for any vertex $x$, 
 the intersection of the ball of radius $r$
centered at $x$ and $C(x,y)$
contains at most $P(r)$ points. 
\end{lem}

{\it Proof of Proposition \ref{prop.raag}}. 
The Cayley graph of $\Gamma$ w.r.t $S$ is the $1$-skeleton 
of the Salvetti complex, $C(L)$,  associated to $L$,
which is a CAT(0) cube complex
(See for example, \cite{NR}.)
The dimension of $C$ is $J$.

Due to Lemma \ref{lemma.CR}, Lemma \ref{lem-mi} applies to the complex $C(L)$ for $G$, 
and we conclude that 
$\Gamma$ satisfies the coarse median inequality of rank $(J+1)$.
\qed

 {\bf Coarse median terminology. }Inequality  (\ref{median-structure}) was termed the spherical coarse median inequality, and the parameter $d_2+1$ associated with it the rank of the inequality.  
To explain this terminology, note first that by Lemma \ref{lem-mi}, a median space satisfies the spherical coarse median inequality. However, by Proposition \ref{Coarse Med Hyp} this inequality is satisfied by all hyperbolic groups w.r.t. word metrics, and these typically do not satisfy the median property, but only a coarser property, which was termed almost median by Bader-Furman \cite{BF17}.  The terminology regarding coarser notions of the median property is however not settled. We note that Bowditch \cite[\S 2]{Bo13} initiated the definition of coarse median spaces and a notion of rank for them. Lafforgue \cite[Def. 2.2]{La00} defines the notion of $H_\delta$-spaces, which was subsequently reformulated by Chatterji-Ruane \cite[Prop. 1.7]{CR05}, and includes also a notion of rank.  Chatterji-Drutu-Haglund \cite[Def. 2.16]{CDH} define the notion of approximate $\delta$-medians and Nica \cite[Def. 4.2]{Ni24} defines the notion of rough median and roughly modular space, similar to those introduced by Lafforgue and Chatterji-Ruane, and associated with a notion of rank. Finally, we mention also the notion of centroids in groups, which are also  relevant to the discussion of coarse properties, and refer 
to \cite{Sa15} for a discussion. It would be desirable to explicate the implications between these notions, and to determine the exact assumptions under which the spherical coarse median inequality (\ref{median-structure}) holds.

\section{Rationality of growth function and  almost exact polynomial-exponential growth}\label{sec.example}

To apply Theorem \ref{thm-HL-d>0} to various groups (in particular, to RAAGs) we need to establish 
almost-exact growth for them. In the present section we discuss a condition of on a finite generating set $S$ of a finitely generated group $G$ which 
gives almost exact polynomial-exponential growth.

\subsection{Rationality of a growth function }

Let $(\Gamma,\abs{\cdot}_S)$ be the Cayley graph of an infinite  group $\Gamma$ w.r.t.  a finite symmetric generating set $S$.
Let $S_n$ be the sphere of radius $n$.
{\em The growth function} is defined as
$$f_{\Gamma,S}(x)=\sum_{n \ge 0} |S_n|x^n$$

A language is {\em regular} if it is the language accepted by 
a finite state automaton whose edges are labeled by elements 
in $S$. See for example, \cite{E}.
Suppose there is a regular language $\mathcal L$ on $S$ that has a bijection 
to the set of elements in $\Gamma$.
If $\mathcal L$ consists of geodesics in $(\Gamma, \abs{\cdot}_S)$, then 
we say $\mathcal L$ is a {\em regular normal form} for $(\Gamma,\abs{\cdot}_S)$.

We quote two relevant facts as a proposition. 
For the first fact, see for example \cite{E}. 
For the second fact, see for example  \cite[\S 3.1, Prop. 3.1.4]{Ca11}, which refers to \cite[Theorem IV.9, Theorem V.3]{FS09}.

\begin{prop}[Rationality and  growth]\label{rational}
Let $\Gamma$ be an infinite  group with a finite symmetric generating set $S$.
Let $(\Gamma,\abs{\cdot}_S)$ be its Cayley graph. 
\begin{enumerate}
\item
Suppose there is a regular normal form for $(\Gamma,\abs{\cdot}_S)$. 
Then $f_{\Gamma,S}$ is rational.  

\item
If in addition $\Gamma$ is of exponential growth, then $\Gamma$
has almost exact polynomial-exponential growth with 
parameters $(q,d)$ for some $q,d$, w.r.t. $\abs{\cdot}_S$.

\end{enumerate}
\end{prop}

Regarding the growth parameters, we note the following. 
 Let $M$ be the matrix associated to the finite state automaton
for the language $\mathcal L$. Then $q$ is equal to  the largest real root of $M$.

We state the following instance of Theorem \ref{main}.
\begin{thm}\label{main2}
Let $\Gamma$ a finitely generated group with exponential growth, $S$ a finite symmetric generating set, and $(\Gamma,\abs{\cdot}_S)$ the associated Cayley graph.
Assume that $\Gamma$  has the rapid decay property w.r.t. $\abs{\cdot}_S$ with polynomial parameter $\bf b$, and in addition has a 
rational growth function. 
Then the Hardy-Littlewood operator satisfies the weak-type  $\cL\left(\log \cL\right)^{2{\bf b}}(\Gamma)$. 
\end{thm}

\proof
Assuming $\Gamma$ has the rapid decay property w.r.t. $\abs{\cdot}_S$ with polynomial parameter $\bf b$, it satisfies
rapid decay of spherical shell correlations for $\bf b$, 
by Proposition \ref{RD-prop}.
Assuming also that $(\Gamma,\abs{\cdot}_S)$ has a rational growth function, 
it has almost exact polynomial-exponential growth by Proposition \ref{rational}.
The result then follows from Theorem \ref{main} applies to $(\Gamma,\abs{\cdot}_S)$.
\qed

\section{Examples of groups satisfying the weak-type maximal inequality}

In this section we give examples of infinite groups that satisfy Hardy-Littlewood 
maximal inequality. First, as a consequence of Theorem \ref{thm-HL-d>0}, we have : 
 
\begin{thm}[RAAG]\label{thm.raag}
Let $\Gamma$ be a RAAG of dimension $J$ with the standard generating set $S$. Suppose that $\Gamma$ is not Abelian.  

 Then the Hardy-Littlewood operator satisfies the weak-type  $\cL\left(\log \cL\right)^{\bf c}(\Gamma)$ maximal inequality, with ${\bf c}=2J+d$,
 where $d$ is the polynomial parameter in the almost exact polynomial-exponential growth.

\end{thm}
\proof
By Proposition \ref{prop.raag}, $(\Gamma,\abs{\cdot}_S)$ satisfies the spherical coarse median inequality of rank $(J+1)$, namely with parameter $d_2=J$. 

It is known that a RAAG has a rational growth function w.r.t. the standard generators $S$, see
\cite[Example 1]{LMW}.

Since $\Gamma$ is not Abelian, it has exponential growth. 
Therefore $(\Gamma,\abs{\cdot}_S)$ has almost exact polynomial-exponential growth by Proposition \ref{rational}.
Now Theorem \ref{thm-HL-d>0} applies to $(\Gamma,\abs{\cdot}_S)$. 
\qed

Note that $d$ is determined by the growth function, and hence $\bf c$ in this case
is an explicit function of the given RAAG data of $\Gamma$.

We now turn to give more examples for groups satisfying the maximal inequality based on Theorem \ref{main2}. Thus the groups we consider satisfy RD, 
and have a rational growth function, which implies almost exact polynomial-exponential 
growth.

If a finitely generated infinite group $\Gamma$ has RD  for some 
length function, then it has RD for any word metric (see \cite{Ch17}), so we simply say $\Gamma$ has RD.

\begin{thm}[More examples for HL maximal inequality]\label{thm.example}
Let $(\Gamma,\abs{\cdot}_S)$ be one of the following groups:
\begin{enumerate}
\item a geometrically finite hyperbolic group 
 (a subgroup  of $Isom(\Bbb H^n)$ by definition)
and a suitable finite generating set $S$ (explained below), 
\item
a Coxeter group of exponential growth 
with the standard generating set $S$,
\item a braid group of exponential growth with 
a suitable standard generating set $S$ (explained below), 
\item an Artin group of extra-large type with 
the standard generating set $S$.

\end{enumerate}

Then $(\Gamma,\abs{\cdot}_S)$ has almost exact polynomial-exponential growth, and 
satisfies the rapid decay property with some parameter ${\bf b}$.  
Therefore, the Hardy-Littlewood operator satisfies the weak-type  $\cL\left(\log \cL\right)^{2{\bf b}}(\Gamma)$ maximal inequality.
\end{thm}

\proof
We check the conditions of Theorem \ref{main2} for each case. 
\\
(1).
By \cite[Theorem 4.3]{NS},
$\Gamma$ has a generating set $S$ with the property called 
``the falsification by fellow traveler property,''
for which the growth function if rational. 
$\Gamma$ has exponential growth, so that $(\Gamma,\abs{\cdot}_S)$
has almost exact polynomial-exponential growth. 

$\Gamma$ satisfies the rapid decay property, as follows from the fact that 
a finitely generated group hyperbolic relative to a family of 
subgroups satisfying RD, itself satisfies RD, see \cite{DS}.
Now any finitely generated group of polynomial growth clearly satisfies RD, and 
the group $\Gamma$ is hyperbolic relative to finitely generated abelian groups, 
so that $\Gamma$ satisfies  RD.

(2). A Coxeter group with the standard generating set
has a rational growth function, \cite[Corollary 1.29]{S}.
It satisfies RD by \cite{CR05}. 

(3) A braid group has a generating set $S$ 
whose growth function is rational, see \cite{CM}.
A braid group satisfies  RD since any mapping class group satisfies RD, see \cite{BM}.

(4) Recall that an Artin group $\Gamma$ is defined for any given  finite simplicial graph $U$ such that 
every edge between vertices $a,b \in V(U)$ is labelled by an integer $m_{ab}\ge 2$.
The Artin group associated to $U$ is the group $A_U$ with the 
set of generators $\{a| a \in V(U)\}$ subject to the relations
$aba \cdots = bab \cdots$ where both sides have word length equal to $m_{ab}$.
$A_U$ is of {\em large type} if all $m_{ab}$ is at least $3$, and of {\em extra-large type}
if they are all at least $4$. 

An Artin group $\Gamma$ of large type has a regular normal form, so that its growth 
function is rational, w.r.t. the standard generating set by \cite[Theorem 4.1]{HR}.
In fact it is shown to be ``shortlex automatic'', which implies that it 
is ``weakly geodesically automatic'', \cite[Corollary 2.5.2]{Ebook}, which 
gives a regular normal form, \cite[Theorem 2.5.1]{Ebook}.

Finally, an Artin group $\Gamma$ of extra-large type satisfies RD, see \cite{CHR2}. 
\qed
\subsection{Other groups}

We turn to briefly discuss some other examples of groups of interest.

Let $\Gamma$ be the mapping class group of a compact 
orientable surface. 
Then $\Gamma$ satisfies RD, see \cite{BM}. 
Indeed, it is prove that $\Gamma$ satisfies RD w.r.t. the metric induced from the ``marking complex'', 
which is quasi-isometric to a word metric. 
It does not seem to be known if $\Gamma$ has almost exact  polynomial-exponential growth 
w.r.t. the metric in question or any other word metrics. 

\begin{question}[Mapping class groups]
Does a mapping class group have almost exact polynomial-exponential growth
w.r.t. a word metric or w.r.t. the metric in \cite{BM}?
\end{question}

Let $X$ be a finite dimensional CAT(0) cube complex.
Let $\Gamma$ be a group that acts on $X$ properly, cellularly and cocompactly.
Then it is known that $\gamma$ satisfies RD, see \cite[Theorem 0.4]{CR05}.
We mention the following example, which
is not a hyperbolic group.
\begin{example}
\begin{enumerate}
 \item
 Let $\Gamma$ be a free by cyclic group of the form $F_2 \rtimes \Z$, then 
$\Gamma$  acts on a CAT(0) square complex, freely, properly,
 and cocompactly, \cite[Theorem 2.1]{BK}.

\end{enumerate}
\end{example}
We remark that any finitely presented subgroup of $\Gamma$ described in the foregoing example acts on a CAT(0) square complex, freely, properly,
 and cocompactly. (This is due to \cite{BH}, see \cite{BK}.)

It would be interesting to know if such a group $\Gamma$ has a word metric or a length function 
for which $\Gamma$ has almost exact polynomial-exponential growth.
We refer to \cite{NR, Sw} for the rationality of a certain language associated to 
$G$.

Finally, let us note that a solvable Baumslag-Solitar group with the standard generators has a rational growth function, see \cite{CEG}, but it does not satisfy RD. It is not known whether it satisfies radial rapid decay with respect to some length function. 

\section{Products}
In the present section we first discuss product groups, and then bring up some compelling problems about optimality of the maximal inequality, which 
are open even in the simple set-up to which we now turn. . 
\subsection{Products}
For $(\Gamma_1,\abs{\cdot}_{S_1}), (\Gamma_2, \abs{\cdot}_{S_2})$ the Cayley graphs of two finitely generated groups, we consider their $\ell_1$-product, namely the group $\Gamma_1 \times \Gamma_2$, taken with the word metric associated with the set of generators is $\{(s,1)|s\in S_1\} \cup \{(1,s)|s \in S_2\}$, 
which we write as $S_1 \cup S_2$. 

We begin by recalling an important stability property of RD under graph-products. 
and refer to \cite{CHR} for a definition of general graph-products. 
\begin{thm}{\cite{CHR} 
The RD property is preserved under graph products, and in particular it is 
preserved under direct products and under free products.}

\end{thm}

It is clear from the definition that the growth function of  the direct product $(\Gamma_1 \times \Gamma_2, S_1\cup S_2)$ is the product of $f_{\Gamma_1,S_1}$ and $f_{\Gamma_2,S_2}$.
Therefore, if both $f_{\Gamma_1,S_1}$ and $f_{\Gamma_2,S_2}$ are rational, then so is $f_{(\Gamma_1 \times \Gamma_2, S_1 \cup S_2)}$.

We state an immediate consequence. 
\begin{cor}[Product]\label{cor.product}
If both $(\Gamma_1,\abs{\cdot}_{S_1})$ and $ (\Gamma_2, \abs{\cdot}_{S_2})$ are of exponential growth, satisfy RD, and have rational growth functions, 
so does their $\ell_1$-product. Hence the Hardy-Littlewood operator 
for $\Gamma_1 \times \Gamma_2$ w.r.t. the metric defined by the generating set $ S_1 \cup S_2$ satisfies the weak-type  $\cL\left(\log \cL\right)^{\bf c}(\Gamma_1\times \Gamma_2)$ maximal inequality for some $\bf c$. 

\end{cor}
\proof
Apply Theorem \ref{main2} to the product.
\qed

For example, we can take non-elementary hyperbolic groups, a non-abelian RAAG 
or any of groups in Theorem \ref{thm.example} as choices for $(\Gamma_1, \abs{\cdot}_{S_1})$ and $(\Gamma_2,\abs{\cdot}_{S_2})$.

\subsection{Products of hyperbolic groups}
By  Corollary \ref{cor.product}, the Hardy-Littlewood maximal operator of the product of two non-elementary 
hyperbolic groups satisfies the weak-type  $\cL\left(\log \cL\right)^{\bf c}$ maximal inequality for some $\bf c$.
Let us address the issue of estimating the constant $\bf c$ in this case.

\begin{lem}[growth of products]\label{lem.growth.product}
Suppose $\Gamma_i$ be arbitrary infinite groups finitely generated groups, with finite symmetric generating sets $S_i$, $i=1,2$. Assume that the spheres $S_i^n$ defined by the length functions $\abs{\cdot}_{S_i}$ have almost exact exponential growth with parameter $q_i > 1$, with  $q_1 \le q_2$. 

Then, $\Gamma_1 \times \Gamma_2$ with the $\ell^1$-product metric associated with the generating set $S_1 \cup S_2$ has almost exact polynomial-exponential growth
of parameters $(d=1,q_1)$ if $q_1=q_2$; and almost exact exponential growth of
parameter $q_2$ if $q_1 < q_2$.

\end{lem}
\proof
Let $S_i^n$ denote the sphere of radius $n$ in the Cayley graph of $(\Gamma_i,\abs{\cdot}_{S_i})$,
where $S_i^0=\{1\}$. (We write $S_i^n$ rather than $(S_i)_n$ in here.)

Let $(S_1 \cup S_2)^n$ denote the sphere of radius $n$ for $(\Gamma_1 \times \Gamma_2, \abs{\cdot}_{S_1 \cup S_2})$. We have
$$(S_1 \cup S_2)^n = \sqcup_{i=0}^n S_1 ^{n-i} S_2^i$$
By assumption, there are constants $C_1,C_2 >0$ such that 
for all $n \ge 0$, we have 
$$\frac{q_i^n}{C_i} \le |S_i^n| \le C_i q_i^n$$
This implies
$$\frac{1}{C_1C_2} \sum _{i=0}^n q_1^{n-i}q_2^i \le |(S_1 \cup S_2)^n| \le
C_1C_2 \sum _{i=0}^n q_1^{n-i}q_2^i$$

(i) Suppose $q_1=q_2$. 
Then 
$$\frac{1}{C_1C_2} (n+1) q_1^{n} \le |(S_1 \cup S_2)^n| \le
C_1C_2 (n+1)q_1^{n},$$
which implies that it has almost exact polynomial-exponential growth of 
parameter $(1,q_1)$.

(ii) Suppose $q_1 < q_2$.
Then, there is a constant $C\ge 1$ such that for all $n \ge 0$, we have 
$$\frac{q_2^n}{C} \le \sum _{i=0}^n q_1^{n-i}q_2^i = \frac{q_2^{n+1}-q_1^{n+1}}{q_2-q_1} \le C q_2^n,$$
which implies almost exact exponential growth of parameter $q_2$. 
\qed

\begin{lem}\label{lem.coarse.product}
Let $\Gamma_i$ be non-elementary hyperbolic groups with finite symmetric generating sets $S_i$. Assume that spheres $S_i^n$, $i=1,2$ have  
almost exact exponential growth with parameters $1< q_1 \le q_2$, respectively. 
Then, $\Gamma_1 \times \Gamma_2$ with the generating set $S_1 \cup S_2$ satisfies rapid decay of 
spherical shell correlations, with polynomial parameter ${\bf b}= 3/2$ when $q_1=q_2$, and with parameter ${\bf b}=2$ when $q_1<q_2$.

\end{lem}
\proof
By a well-known result of Coornaert \cite{Co93} for spheres in a word-hyperbolic,   there exist constants  $C_1,C_2$ such that for any $n \ge0$, 
$$ \frac{q_i^n}{C_i} \le |S_i^n| \le C_i q_i^n.$$
Suppose that the Cayley graphs for $(\Gamma_i,\abs{\cdot}_{S_i})$ are $\delta$-hyperbolic. 
Put $B_i=q_i^{10\delta}$.
Let $C =\max\{C_i\}, B=\max\{B_i\}$, and $q=\max\{q_i\}$. 

Let $E,F \subset \Gamma_1\times \Gamma_2$ be finite subsets.

We will show 
\begin{equation}\label{eq:split levels 3} |\{ (x,y) \in E \times F:
d_{\abs{\cdot}_{S_1\cup S_2}}(x,y) = r \}| \le 2B^2 C^2
\sqrt{|E|} \sqrt{|F|}  r^2 \sqrt{q^r}.
\end{equation}

Note that this inequality gives the desired conclusion. Indeed, if $q_1=q_2=q$, then
by Lemma \ref{lem.growth.product}, for any $r>0$
 $$ \frac{r q^r}{C_1} \le (S_1 \cup S_2)^r \le C_1r q^r,$$ where $C_1$ is a constant.
 The first inequality implies 
$$ r^2 \sqrt{q^r} \le  \sqrt{C_1} r^{3/2} \sqrt{(S_1 \cup S_2) ^r }.$$
Substituting this into (\ref{eq:split levels 3}), we get the desired estimate with ${\bf b}=3/2$.

Also, if $q_1 < q_2 =q$, then 
by Lemma \ref{lem.growth.product}, 
 $$  \frac{q^r}{C_2} \le (S_1 \cup S_2)^r \le C_2  q^r,$$ where $C_2$ is a constant. Hence
$$ r^2 \sqrt{q^r} \le  \sqrt{C_2} r^2 \sqrt{(S_1 \cup S_2) ^r }, $$
which gives the desired estimate with ${\bf b}=2$.

We turn to proving (\ref{eq:split levels 3}), and refer to the proof of Proposition \ref{Coarse Med Hyp} and Fig. 1. 

We denote the distance $\abs{\cdot}_{S_i}$ on $\Gamma_i$ by  $d^{(i)}$, and the distance on 
the product group $\Gamma_1\times \Gamma_2$ by $d^{(p)}$. 

As before, set 
$$E_n=E \cap (S_1 \cup S_2)^n, F_n= F \cap (S_1 \cup S_2)^n.$$

Given integers $i,j,r\ge0$ satisfying the triangle inequalities, set $$m=\frac{j+r-i}{2}$$
We have $ 0 \le m \le r$. 

Let us rewrite the main expression as :
\begin{equation}\label{eq:split levels} |\{ (x,y) \in E \times F:
d^{(p)}(x,y)  = r \}| 
= \sum_{m=0}^r \sum_{\substack{i,j \in \N\cup\{0\}\\
i = j+r-2m}} \left| \{ (x,y) \in E_j \times F_{i}: d^{(p)}(x,y)  = r
\}\right|.\end{equation}

We first prove
\begin{equation}\label{eq.product}
\left| \{ (x,y) \in E_j \times F_{i}: d^{(p)}(x,y)  = r
\}\right| \leq B^2C^2  r^2 \min\left\{ q^{r-m} |E_j|, q^m |F_{i}|
\right\}.\end{equation}

Write $x=(x_1,x_2) \in E_j$ and 
$y=(y_1,y_2)\in F_i$.  Set $j_1=d^{(1)}(1,x_1), j_2=d^{(2)}(1,x_2)$ and $i_1=d^{(1)}(1,y_1), i_2=d^{(2)}(1,y_2)$. 

We have 
$$|\{ (x,y) \in E_j \times F_{i}: d^{(p)}(x,y) = r
\} | $$
$$= \sum_{0\le r_1 \le r,}
\sum_{  0 \le m_1 \le r_1} |\{ (x,y) \in E_j \times F_{i}: d^{(p)}(x,y) = r, d^{(1)}(x_1,y_1)=r_1, 
i_1= r_1+ j_1-2m_1
\}|$$

Set $r_2=r-r_1$ and $m_2=m-m_1$. 

\begin{figure}[h]
\begin{center}
\includegraphics[scale=0.45, angle=90]{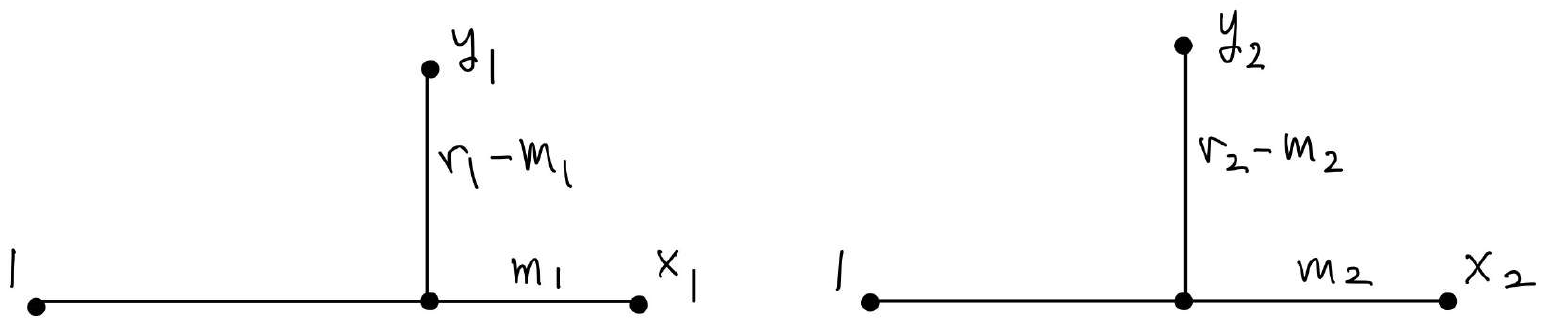}
   \end{center}
\caption{Schematic pictures for the $\delta$-thin triangles }
  \label{fig.product}
\end{figure}

Now, fix a point $x=(x_1,x_2) \in E_j$. Given $r_1,m_1$ as above, 
we want to count the possibilities for $y=(y_1,y_2)\in F_i$
with $$d^{(p)}(x,y)=r,  d^{(1)}(x_1,y_1)=r_1, 
d^{(1)}(1,y_1)=i_1= r_1+ j_1-2m_1$$ 
By the $\delta$-hyperbolicity of the Cayley graph of $(\Gamma_1,\abs{\cdot}_{S_1})$, the number of possibilities
for such $y_1$ is at most $CBq_1^{r_1-m_1}$.
Also, by the $\delta$-hyperbolicity of the Cayley graph of $(\Gamma_2,\abs{\cdot}_{S_2})$,  the number of possibilities
for such $y_2$ is at most $CBq_2^{r_2-m_2}$,  
see Figure \ref{fig.product}.

Thus the possibilities for such $y=(y_1,y_2)$ is at most
$$C^2B^2q_1^{r_1-m_1} q_2^{r_2-m_2} \le C^2B^2q^{r_1+r_2-m_1-m_2}=C^2B^2q^{r-m}$$
Therefore 
$$\left| \{ (x,y) \in E_j \times F_{i}: d^{(p)}(x,y) = r, d^{(1)}(x_1,y_1)=r_1, 
i_1= r_1+ j_1-2m_1
\}\right| 
$$
$$
\leq B^2C^2    q^{r-m} |E_j|.
$$

It follows, since $0 \le m_1 \le r_1$ and $0 \le r_1 \le r$, 

$$\left| \{ (x,y) \in E_j \times F_{i}: d^{(p)}(x,y) = r
\}\right| \leq B^2C^2  r^2 q^{r-m} |E_j|
$$
Similarly, we have
$$\left| \{ (x,y) \in E_j \times F_{i}: d^{(p)}(x,y) = r
\}\right| \leq B^2C^2  r^2 q^{m} |F_i|
$$
Combining those two inequalities, we proved (\ref{eq.product}).

Inserting (\ref{eq.product}) into (\ref{eq:split levels}), we get
\begin{equation}\label{eq:split levels 2} |\{ (x,y) \in E \times F:
d^{(p)}(x,y) = r \}| \le B^2 C^2 r^2 \sum_{m=0}^r \sum_{\substack{i,j \in \N\cup\{0\}\\
i = j+r-2m}}
\min\left\{ q^{r-m} |E_j|, q^m |F_{i}|
\right\}.
\end{equation}

To bound the latter expression, we use the argument used in the proof of Proposition \ref{convolution-ineq} which as noted above generalizes the method developed in \cite[Proof of Lemma 5.1]{NT09} for free groups), and conclude :  
\begin{equation}\label{eq:tree goal} \sum_{m=0}^r \sum_{\substack{i,j \in \N\cup\{0\}\\
i = j+r-2m}} \min\left\{ q^{r-m} |E_j|, q^m |F_{i}| \right\} \le
2\sqrt{|E|} \sqrt{|F|} \sqrt{q^r}.\end{equation}

Inserting  (\ref{eq:tree goal})   in (\ref{eq:split levels 2}),
we get (\ref{eq:split levels 3}). 
\qed

\begin{thm}\label{thm.HLM.product}
Let $\Gamma_i$ be non-elementary hyperbolic groups, $S_i$ finite symmetric generating sets, $i=1,2$. Assume $(\Gamma_1,\abs{\cdot}_{S_1})$ has almost exact exponential growth of parameter $q_1>1$ and 
$(\Gamma_2,\abs{\cdot}_{S_2})$ has almost exact exponential growth of parameter $q_2>1$ with $q_1 \le q_2$. 
Then,
\begin{enumerate}
\item
if $q_1=q_2$, then $(\Gamma_1 \times \Gamma_2, S_1 \cup S_2)$ has almost exact polynomial-exponential growth parameters $(1,q_1)$ (namely with polynomial parameter $d=1$). 
The Hardy-Littlewood operator satisfies the weak-type  $\cL\left(\log \cL\right)^{\bf c}(\Gamma_1 \times \Gamma_2)$ maximal inequality, with ${\bf c}=3$.

\item
If $q_1 < q_2$, then $(\Gamma_1 \times \Gamma_2, S_1 \cup S_2)$ has almost exact exponential growth
of parameter $q_2$, and it satisfies the coarse median inequality with parameter $d_2=2$.
The Hardy-Littlewood operator satisfies the weak-type  $\cL\left(\log \cL\right)^{\bf c}(\Gamma_1 \times \Gamma_2)$ maximal inequality, with ${\bf c}=2(d_2+0)=4$. 

\end{enumerate}

\end{thm}
\proof
By Lemma \ref{lem.coarse.product} and Lemma \ref{thm.HLM.product}, 
Theorem \ref{main} applies to $(\Gamma_1 \times \Gamma_2, S_1 \cup S_2)$
with ${\bf c}=2{\bf b}$, 
and we get the conclusion with ${\bf b} = 3/2$ in the first case,
and ${\bf b} = 2$ in the second case. 

\qed

\subsection{The problem of optimality}

 Consider length functions on an infinite finitely generated group $\Gamma$, with almost exact polynomial-exponential growth of spherical shells, and note the following. 
 
For word-metrics associated with a finite symmetric generating set,  in \cite{Ni10} it is shown that when the rapid decay property holds necessarily the polynomial parameter satisfies ${\bf b}\ge \frac12$. 
It follows that when a maximal inequality in $\cL\left(\log \cL\right)^{\bf c}(\Gamma)$ is proved for the Hardy-Littlewood operator using the {\it rapid decay property}, necessarily  ${\bf c}\ge 1$. 

It follows from Proposition \ref{convolution-ineq} that for integer-valued length functions,  when the spherical coarse median inequality with polynomial parameter $d_2=0$ holds together with almost exact {\it exponential} growth, then rapid decay of spherical correlations holds with polynomial parameter ${\bf b}=0$. 
In Proposition \ref{Coarse Med Hyp} it was established that word metrics on hyperbolic  groups satisfy the coarse median inequality with parameter $d_2=0$. As stated in Theorem \ref{thm.hyp},  in that case the weak-type $(1,1)$-maximal inequality in $\ell^1(\Gamma)$ holds for the Hardy-Littlewood operator, which is the  {\it best-possible} result. 
 
However, {\it for all other examples} presented in the present paper for which the maximal inequality in $\cL\left(\log \cL\right)^{\bf c}(\Gamma)$ was established for some length function, the optimal value of $\bf c$ is an open problem. It may be the case that the best-possible result ${\bf c}=0$ holds in all of these examples.

What then should one expect about the optimal value of $\bf c$ for which the Hardy-Littlewood maximal inequality is valid ? First note that for a connected semisimple Lie group $G$ with finite center and the Riemannian balls $\tilde{B}_r$, the optimal value ${\bf c}=0$ was established in Str\"omberg's theorem \cite{Str81}, as noted in the introduction. Perhaps this fact can be taken as an indication that for the balls $\tilde{B}_r\cap \Gamma$ on any lattice $\Gamma \subset G$ the optimal result holds. It may also be possible, in principle, that the result may be different for uniform lattices and for non-uniform lattices, for the case of Riemannian ball averages or for other length functions. 

For ball averages with respect to natural distances on Bruhat-Tits buildings associated with a linear algebraic group, the optimality problem for the Hardy-Littlewood operator is completely open. Currently the only exceptions are the case of regular and bi-regular trees. 
  Optimality for the Hardy-Littlewood operator on lattices in the associated algebraic group, with respect to the restricted metrics from $G$ or word-metrics, is likewise completely open, again other than the case when their associated building is a regular or bi-regular tree.


\begin{thebibliography}{BKK11}

\bibitem[An87]{An87} J.-P. Anker, \textit{La forme exacte de l’estimation fondamentale de Harish-Chandra}, C. R. Acad. Sci. Paris Ser. I Math. 305 (1987), 371- 374. 



\bibitem[ADY96]{ADY96} J-P. Anker, E. Damek and Ch. Yacoub, \textit{Spherical analysis on harmonic $AN$ group}. Annali Scuola Nor. Sup. di Pisa, vol. 23, pp. 643-679 (1996). 

\bibitem[BF17]{BF17} U. Bader , and A. Furman, \textit{Some ergodic properties of metrics on hyperbolic groups}.  Math. ArXiv:1707.02020, August 2017.

\bibitem[BM]{BM}
J. A. Behrstock and  Y. N. Minsky, \textit{Centroids and the rapid decay property in mapping class groups}. J. Lond. Math. Soc. (2) 84 (2011), no. 3, 765-–784. 



\bibitem[BR]{BR}
J. Berstel and Ch. Reutenauer, \textit{Another proof of Soittola's theorem}. Theoret. Comput. Sci. 393 (2008), no. 1--3, 196–-203. 

\bibitem[BS93]{BS93} C. Bishop and T. Steger, \textit{representation-theoretic rigidity 
in $PSL_2(\RR)$}. Acta Math. vol. 170 (1993) pp. 121-149. 


\bibitem[B24]{B24} K. Boucher, \textit{Spectral gaps, critical exponents, and representations of negatively curved groups.} Math. ArXiv 2401. 16962, Jan. 2024

\bibitem[Bo13]{Bo13} B. H. Bowditch, \textit{Coarse median spaces and groups}. Pacific Journal of Mathematics, vol. 261, pp. 53-93, 2013. 


\bibitem[BLP17]{BLP17} A. Boyer, G. Link, and C. Pittet, \textit{ Ergodic properties of boundary representations.} Ergod. Th. and Dynam. Sys., 39:2023–2047, 2017.

\bibitem[Bo13]{Bo13} B. H. Bowditch, \textit{Coarse median spaces and groups}. Pacific Journal of Mathematics, vol. 261, pp. 53-93, 2013. 

\bibitem[BLP22]{BLP22} A. Boyer, A. P. Lobos and Ch. Pittet, \textit{Radial rapid decay does not imply rapid decay}. Annales de l'Institut Fourier, vol. 73, pp. 1421-1452, 2023. 

\bibitem[BH]{BH}
M.R. Bridson,  A. Haefliger, \underline{Metric Spaces of Non-Positive Curvature}. 
Grundlehren Math.  319
Springer-Verlag,  1999.  
%

\bibitem[Br05]{Br05} E. Breuillard, \textit{Geometry of locally compact groups of polynomial growth and shape of large balls.} Groups, geometry and Dynamics, vol. 8
pp. 639-732, 2014. 

\bibitem[BK]{BK}
Button, J. O. ; Kropholler, R. P., \textit{Nonhyperbolic free-by-cyclic and one-relator groups}. 
New York J. Math. 22 (2016), pp. 755–774. 


\bibitem[C53]{C53} Calderon, A., \textit{ A general ergodic theorem}. Ann. of Math. 57  (1953), pp. 182-191. 


 \bibitem[Ca11]{Ca11} Calegari, D., \textit{The ergodic theory of  hyperbolic groups}. Contemp. Math. 597 (2013), pp. 15--52. 
 


\bibitem[CM]{CM}
Charney, Ruth; Meier, John, \textit{ The language of geodesics for Garside groups}. Math. Z. 248 (2004), no. 3, 495–509.



\bibitem[Ch17]{Ch17} I. Chatterji, \textit{Introduction to the Rapid Decay property}.  Contemporary Mathematics
Vol. 691, 2017

\bibitem[CDH]{CDH} I. Chatterji, C. Drutu and F. Haglund, \textit{Kazhdan and Haagerup properties from the median viewpoint. } Advances in Mathematics, vol. 225, pp. 882-921, 2010. 


\bibitem[CPSC]{CPSC} I. Chatterji, C.  Pittet and L. Saloff Coste, \textit{Conected Lie groups and  Property RD}. Duke Math. Jnl. 
Vol. 138 (2007), pp. 511-536. 



\bibitem[CR05]{CR05}
I. Chatterji, K. Ruane, \textit{Some geometric groups with rapid decay}.
Geom. Funct. Anal. vol. 15 (2005), no. 2, 311–-339. 




\bibitem[CHR]{CHR} L. Ciobanu, D. Holt, Derek F. and S. Rees,  \textit{Rapid decay is preserved by graph products}. J. Topol. Anal. 5 (2013), no. 2, 225–237. 
 
 \bibitem[CHR2]{CHR2} L. Ciobanu, D. Holt, Derek F. and S. Rees, \textit{Rapid decay and Baum-Connes for large type
Artin groups}. Trans. Amer. Math. Soc.368(2016), no.9, 6103--6129.




\bibitem[CEG]{CEG}
D. J. Collins, M.  Edjvet, M. and C. P. Gill, \textit{Growth series for the group $<x,y| x^{-1}yx=y^\ell>$}.
Arch. Math. 62 (1994), no. 1,pp.  1–11. 

\bibitem [Co93]{Co93} M. Coornaert,  \textit{Mesures de Patterson-Sullivan
 sur le bord d'un espace hyperbolique au sens de Gromov}.
Pacific J. Math., vol. 159  no. 2, pp.  241-270 (1993).

\bibitem[Co78]{Co78} M. Cowling, \textit{The Kunze–Stein phenomenon}, Ann. Math. 107 (1978), 209-234.

\bibitem[CMS]{CMS} M. G. Cowling, S. Meda, and A. G. Setti, \textit{ A weak-type $(1,1)$ estimate for a maximal operator on a group of isometries of homogeneous trees}. 
Colloquium Mathematicae,  vol. 118 (2010) pp. 223-232. 

\bibitem[CS96]{CS96} M. Cowling and T. Steger, \textit{Irreducibility and restrictions of unitary representations to lattices}. J. Reine. Angew. Math (Crelle Journal), vol. 420, pp. 85-98 (1996). 

\bibitem[DRS93]{DRS93} W. Duke, Z. Rudnick, P. Sarnak, \textit{Density of integer points on affine  homogeneous varieties}.
Duke Math. J. 71 (1993), no. 1, 143–179.

\bibitem[DS]{DS}
C. Drutu, and M. Sapir, \textit{Relatively hyperbolic groups with rapid decay property}.
Int. Math. Res. Not. 2005, no. 19, 1181–1194. 



\bibitem[E]{E}
 D.B. A. Epstein, A. R. Iano-Fletcher, Uri  Zwick, \textit{Growth functions and automatic groups}. 
Experiment. Math. 5 (1996), no. 4, 297-–315.

\bibitem[CEHPT]{Ebook}
 D.B.A. Epstein, J.W. Cannon, D.F. Holt, S.V.F Levy,; M.S. Paterson, W.P.  Thurston,  \underline{\it Word Processing in Groups}. Jones and Bartlett Publishers, 1992. 

\bibitem[EM93]{EM93} A. Eskin and C. McMullen, \textit{Mixing, counting and equidistribution in Lie groups}. Duke Math. J. 71(1): 181-209 (July 1993)

\bibitem[FS09]{FS09} Ph.  Flajolet, and R. Sedgewick, \underline{\it Analytic Combinatorics}. Cambridge University Press, Cambridge, 2009. 




%
\bibitem[Fa72]{Fa72} A. N. Fava, \textit{Weak type inequalities for product operators}. Studia Math., vol. XLIL, pp. 271-288 (1972). 

\bibitem[Fa84]{Fa84} A. N. Fava, \textit{Mapping properties of maximal operators}. Studia Math. vol. LXXVIII, pp. 1-6 (1984). 

\bibitem[GV88]{GV88} R. Gangolli and V. S. Varadarajan, \underline{\it Harmonic Analysis of Spherical Functions on} \underline{ \it Real Reductive Groups}. A Series of Modern Survey in Mathematics, vol. 101, 1988, Springer Verlag. 


\bibitem[GN10]{GN10} A. Gorodnik and A. Nevo, 
\underline{\it The Ergodic Theory of Lattice Subgroups},
Annals of Mathematics Studies 172, Princeton University Press, 2010.

\bibitem[GN12a]{gn_counting} A. Gorodnik and A. Nevo, \textit{Counting lattice points}, 
 J. Reine Angew. Math. 663 (2012), pp. 127--176. 

\bibitem[GO07]{GO07} A. Gorodnik and H. Oh, \textit{Orbits of discrete subgroups on a symmetric space and the Furstenberg boundary}. 
Duke Mathematical Journal Vol. 139, Issue 3, (2007) , pp.  483-525. 


\bibitem[GW07]{gw}
A. Gorodnik and B. Weiss, \textit{Distribution of lattice orbits on homogeneous varieties}.
Geom. Funct. Anal. 17 (2007), no. 1, pp. 58--115.




\bibitem[Gu73]{Gu73} Y. Guivarc’h, \textit{Croissance polynomiale et periodes des fonctions harmoniques}. Bull. Soc.
Math. France, 101 (1973) pp. 353-379.




\bibitem[HL30]{HL30} G. H. Hardy and J. E.  Littlewood,  \textit{A maximal 
theorem with function theoretic applications}. Acta Math. \textbf{54} (1930),
pp. 81-116.



\bibitem[HR]{HR}
D. F. Holt and S Rees,  \textit{Artin groups of large type are shortlex automatic with regular geodesics}.
Proc. LMS, vol. 104, no. 3, 2012, 486-512.




\bibitem[J]{J} P. Jolissaint, \textit{Rapidly decreasing functions in reduced $C^{\*}$-algebras of groups}, Trans. Amer. Math. Soc. 317 (1990), no. 1, 167–196.

\bibitem[Kn95]{Kn95} G. Knieper, \textit{On the asymptotic geometry of nonpositively curved manifolds}. Geom. Funct. Analysis. Vol. 7 (1997) pp. 755-782. 





\bibitem[LMSV]{LMSV} M. Levi, S. Meda, F. Santagati and M. Vallarino, \textit{Hardy–Littlewood maximal operators on trees with bounded geometry}.
 Math.arXiv:2308.07128v1, to appear in Trans. Amer. Math. Soc. 

\bibitem[La00]{La00} V. Laﬀorgue, \textit{A proof of property (RD) for cocompact lattices of SL(3, R) and SL(3, C).} 
J. Lie Theory vol. 10, pp.  255–267, 2000.

\bibitem[LMW]{LMW}
J. Loeffler, J. Meier and J.  Worthington, \textit{Graph products and Cannon pairs}. 
Internat. J. Algebra Comput. 12 (2002), no. 6, 747–-754. 

\bibitem[LMR00]{LMR00} A. Lubotzky, S. Mozes, and M.S. Raghunathan, \textit{ The word and Riemannian metrics on lattices of semisimple groups}. Publications Matematiques de l’Institut des Hautes Scientifiques, vol. 91, pp. 5-53, 2000.

\bibitem[LZ23]{LZ23} H-Q Li and J-X Zhu, \textit{Weak-type $(1,1)$ of Riesz transform os some direct product manifolds with exponential volume growth}. Math ArXiv:2304.04335v1, April 2023. 


\bibitem[Mac82]{Mac82} H.D. Macpherson, \textit{Infinite distance-transitive graphs of finite valency}, Combinatorics 2 (1982) 63-69.

\bibitem[Ma91]{Ma91} G. A. Margulis, \underline{Discrete Subgroups of Semisimple Lie Groups}, Ergeb. Math. Grenzgeb. (3) 17, Springer, Berlin, 1991.
%

\bibitem[MPSV]{MPSV} S. Meda, S. Pigola, A. G. Setti, and G. Veronelli, \textit{Hardy-Littleewood maximal operators on certain manifoldwith bounded geometry}. 
Math. ArXiv:2502.13109.v1, February 2025. 


\bibitem[NT09]{NT09} Naor, A. and Tao, T., \textit{Random Martingales and localization of maximal inequalities}. J. Funct. Analysis, 2009

\bibitem[NS]{NS}
W. D. Neumann and M. Shapiro, \textit{Automatic structures, rational growth, and geometrically finite hyperbolic groups}.
Invent. Math. 120 (1995), no. 2, 259-–287. 

\bibitem[N98]{N98} A. Nevo, \textit{Spectral transfer and pointwise ergodic theorems for semisimple Kazhdan groups}.  Math. Research Letters, vol. 5, pp. 305-325. 




\bibitem[N04]{N04} A. Nevo,  \textit{Radial geometric analysis on groups}. Discrete Geometric Analysis, Proceedings of the first JAMS conference, held at Tohoku University, Japan, Dec. 2002. Contemp. Math., vol. 347, pp. 221-244 (2004)


\bibitem[N06]{N06} A. Nevo, \textit{Pointwise ergodic theorems for actions of Lie groups}. Handbook of Dynamical systems, vol. 1B, pp. 871-972 , Eds. A. Katok and B. Hasselblatt, Elsevier,  Amsterdam 2006. 





\bibitem[NR]{NR}
 G. A. Niblo, L. D. Reeves, \textit{The geometry of cube complexes and the complexity of their fundamental groups}.
Topology 37 (1998), no. 3, 621–-633. 




\bibitem[Ni10]{Ni10} B. Nica, \textit{On the degree of rapid decay}.  Proc. Amer. Math. Soc. 138 (2010), no. 7, pp. 2341-2347

\bibitem[Ni17]{Ni17} B. Nica, \textit{On operators norms for hyperbolic groups}. Journal of Topology and Analysis, vol. 9, No. 02, pp. 291-296 (2017)

\bibitem[Ni24]{Ni24} B. Nica, \textit{Norms of spherical averaging operators 
for some geometric group actions.} Math. ArXiv, 2405.08682, May 2024. 





\bibitem[ORRS]{ORRS} S. Ombrosi, I.P. Rivera-Rios, and M. D. Safe, {\textit Fefferman-Stein inequalities for the Hardy-Littlewood maximal function of the infinite rooted $k$-ary tree}. IMRN, vol. 2021, pp. 2736-2762, 2021. 


\bibitem[Pa83]{Pa83} Pansu, P.,  \textit{Croissance des boules et des geodesiques fermes dans les nilvarietes}. Erg. Th. Dyn. Sys. 3 (1983), pp. 415-445.


\bibitem[Pe09]{Pe09} Perrone, M. \textit{Radial rapid decay property for co-compact lattices}. 
Journal of Functional Analysis, vol. 256, pp. 3471-3489 (2009). 



\bibitem[RT91]{RT91} R. Rochberg and M. Taibleson, \textit{Factorization of the Green's operator and weak-type estimates for a random walk on a tree}.
Publ. Mat. vol. 35, pp. 187-207 (1991).  

\bibitem[Sa15]{Sa15} M. Sapir, \textit{The rapid decay property and centroids in groups}. 
Journal of Topology and Analysis, vol. 7, No. 03, pp. 513-541 (2015).

  
  \bibitem[ST19]{ST19} J. Soria and P. Tradecete, {\textit Geometric properties of infinite graphs and the Hardy-Littlewood maximal operator}. Journale d'Analyse Mathematiques, 
vol. 137, pp. 913-937, 2019. 

\bibitem[St93]{St93} E. M. Stein, \underline{\it Harmonic Analysis: Real-Variable Methods, Orthogonality and Oscillatory} \underline{Integrals}. Princetom Mathematical Series, vol. 43, Princeton University Press, 1993. 


\bibitem[S]{S} R. Steinberg, \textit{Endomorphisms of Linear Algebraic Groups}.
  Memoirs of the American Mathematical Society
Volume: 1; 1968.

\bibitem[Str81]{Str81} Str\"omberg, J. L., \textit{Weak-type $L^1$-estimates for maximal functions on non-compact symmetric spaces}. Ann. Math. vol.114, pp. 115-126 (1981). 


\bibitem[Sw]{Sw}
 Swiatkowski, Jacek.
 \textit{Regular path systems and (bi)automatic groups}.
Geom. Dedicata 118 (2006), 23–48.


\bibitem[Va97]{Va97} A. Valette, \textit{On the Haagerup inequality and group acting on $\tilde{A}_n$-buildings}.  Ann. Inst. Fourier (Grenoble) 47 (4)
(1997) 1195-1208.


\bibitem[Wi39]{Wi39} Wiener, N. \textit{The ergodic theorem}. Duke Math. J. \textbf{5} (1939), 1--18. 


\end{thebibliography}
\end{document}